\documentclass{amsart}
\usepackage{amssymb,amscd,amsthm,verbatim}
\input epsf                

\newtheorem{proposition}{Proposition}[section]
\newtheorem{theorem}[proposition]{Theorem}
\newtheorem{corollary}[proposition]{Corollary}
\newtheorem{lemma}[proposition]{Lemma}

\newtheorem{remark}[proposition]{Remark}

\newtheorem{prop}[proposition]{Proposition}

\newcommand{\reals}{\mathbb R}
\newcommand{\K}{\mathbb K\,}

\newcommand{\supp}{\mathrm{supp}}

\newcommand{\A}{{\mathcal A}}
\newcommand{\R}{{\mathcal R}}
\newcommand{\B}{{\mathcal B}}

\newcommand{\C}{{\mathcal C}}
\newcommand{\scrH}{{\mathcal {H}}}
\newcommand{\F}{{\mathcal F}}
\newcommand{\Po}{{\mathcal P}}

\newcommand{\covers}{{\,\,\,\cdot\!\!\!\! >\,\,}}
\newcommand{\covered}{{\,\,<\!\!\!\!\cdot\,\,\,}}

\newcommand{\0}{{\hat{0}}}
\newcommand{\1}{{\hat{1}}}
\newcommand{\set}[1]{{\left\lbrace #1 \right\rbrace}}

\newcommand{\join}{\vee}
\newcommand{\meet}{\wedge}
\newcommand{\Cg}{\mathrm{Cg}}
\newcommand{\Con}{\mathrm{Con}}

\newcommand{\Irr}{\mathrm{Irr}}

\newcommand{\st}{\mathrm{st}}

\newcommand{\pidown}{\pi_\downarrow}
\newcommand{\piup}{\pi^\uparrow}

\newcommand{\ep}{\epsilon}

\newcommand{\br}[1]{\langle #1 \rangle}
\newcommand{\Tr}{\mathrm{Tr}}

\newcommand{\LI}{\mathrm{L}}
\newcommand{\RI}{\mathrm{R}}
\newcommand{\product}{\,\substack{\bullet}\,}
\newcommand{\productS}{\,\,{\substack{\bullet}}_{\substack{S}}\,\,}
\newcommand{\productZ}{\,\,{\substack{\bullet}}_{\substack{Z}}\,\,}
\newcommand{\productY}{\,{\substack{\bullet}}_{\substack{Y}}\,}

\begin{document}
\title[Fans and Hopf algebras]{Lattice congruences, fans and Hopf algebras}

\author{Nathan Reading}
\address{
Mathematics Department\\
       University of Michigan\\
       Ann Arbor, MI 48109-1109\\
USA}
\thanks{The author was partially supported by NSF grant DMS-0202430.}
\email{nreading@umich.edu}
\urladdr{http://www.math.lsa.umich.edu/$\sim$nreading/}

\begin{abstract}
We give a unified explanation of the geometric and algebraic properties of two well-known maps, 
one from permutations to triangulations, and another from permutations to subsets.
Furthermore we give a broad generalization of the maps.
Specifically, for any lattice congruence of the weak order on a Coxeter group we construct a complete fan of convex cones with strong properties 
relative to the corresponding lattice quotient of the weak order.
We show that if a family of lattice congruences on the symmetric groups satisfies certain compatibility conditions then the family defines a sub 
Hopf algebra of the Malvenuto-Reutenauer Hopf algebra of permutations.
Such a sub Hopf algebra has a basis which is described by a type of pattern avoidance.
Applying these results, we build the Malvenuto-Reutenauer algebra as the limit of an infinite sequence of smaller algebras, where the second 
algebra in the sequence is the Hopf algebra of non-commutative symmetric functions.
We also associate both a fan and a Hopf algebra to a set of permutations which appears to be equinumerous with the Baxter permutations.
\end{abstract}

\maketitle

\section{Introduction}
\label{intro}
The results of this paper are motivated by the relationship between the permutohedron, the associahedron and the cube, and the corresponding 
relationship between the Malvenuto-Reutenauer Hopf algebra~\cite{MR} of permutations, the Hopf algebra of planar binary trees~\cite{LR} 
and the Hopf algebra of non-commutative symmetric functions~\cite{Gel et al}.
There is a well-known~\cite{Nonpure II,LRorder,EC1,Tonks} map $\eta$ from permutations to Catalan-objects, which has interesting 
properties with respect to these polytopes and algebras.
More precisely, several maps have been studied, related by natural bijections on permutations such as the inverse map, but for the purposes
of this introduction, we call all of these maps ``the map $\eta$.''
In~\cite{Iterated}, Billera and Sturmfels give a realization of the associahedron and the permutohedron such that the normal fan of the 
permutohedron refines that of the associahedron, and $\eta$ is the inclusion map from maximal normal cones of the permutohedron to 
maximal normal cones of the associahedron.
The descent map, mapping a permutation to its descent set, can be realized as the inclusion map from the maximal normal cones of the 
permutohedron to the maximal normal cones of a combinatorial cube.
This map factors through $\eta$, giving a triangle of maps relating the permutohedron to the cube, via the associahedron.
On the algebraic side, the dual maps to this triangle of maps give an embedding of the Hopf algebra of non-commutative symmetric 
functions as a sub Hopf algebra of the Hopf algebra of planar binary trees, and an embedding of the Hopf algebra of planar binary trees as 
a sub Hopf algebra of the Malvenuto-Reutenauer Hopf algebra~\cite{LR}.

The fact that these maps have such nice properties with respect both to polytopes and to algebras demands a unified explanation.
We provide a unified explanation using lattice congruences.
The key to the explanation is the observation that $\eta$ is a lattice homomorphism from the weak order on $S_n$ to the Tamari 
lattice, and that the descent map is a lattice homomorphism from $S_n$ to a Boolean algebra.
A generalization of this observation about the Tamari lattice is proven in~\cite{Cambrian}, although essentially all the ingredients for proving 
it for the Tamari lattice were previously obtained in~\cite{Nonpure II}.
The fact that the descent map is a lattice homomorphism is due to Le Conte de Poly-Barbut~\cite{Barbut}.

The symbol $W$ denotes a finite Coxeter group equipped with the weak order, and $\F$ is the complete fan defined by 
a corresponding Coxeter arrangement.
The combinatorics of the weak order is closely connected to the geometry of~$\F$.
To generalize this close connection, we introduce fan posets and establish their basic properties. 
A {\em fan poset} $(\F,P)$ is a partial order $P$ on the maximal cones of a complete fan $\F$ in $\reals^d$, with some conditions 
relating the partial order to the structure of the fan.
A complete fan $\F$ defines a cellular sphere $\Delta$.
Not every fan is the normal fan of a polytope, but every complete fan has a dual cellular sphere~$\Gamma$ which plays the role of the polytope.
If $(\F,P)$ is a fan poset, then the Hasse diagram of $P$ is isomorphic as a graph to the 1-skeleton of $\Gamma$.
{\em Facial} intervals of $(\F,P)$ are intervals $I$ in $P$ such that, for some cone $F$ of $\F$, the interval $I$ consists of all of the 
maximal cones of $\F$ containing $F$.
The fan poset $(\F,P)$ is {\em homotopy facial} if all non-facial intervals are contractible and if, for every cone $F$, the facial interval 
corresponding to $F$ is homotopy equivalent to a sphere of dimension $d-2-\dim F$.
If $(\F,P)$ is homotopy facial then in particular the M\"{o}bius function of a non-facial interval is zero and the M\"{o}bius function of a 
facial interval corresponding to a face~$F$ is $(-1)^{d-2-\dim F}$.
The fan poset $(\F,P)$ is {\em atomic-facial} if the facial intervals are exactly the atomic intervals.
The definition of {\em bisimplicial} fan posets is given in Section~\ref{fan posets}.

The main geometric result of this paper is a more general version of the following theorem, in which $W/\Theta$ denotes the quotient mod 
$\Theta$ of the weak order on $W$.
\begin{theorem}
\label{main}
If $W$ is a Coxeter group with associated fan $\F$ then for any lattice congruence $\Theta$ on the weak order on $W$ there is a fan 
$\F_\Theta$, refined by $\F$, such that $(\F_\Theta,W/\Theta)$ is a fan lattice.
Furthermore, $(\F_\Theta,W/\Theta)$ is homotopy facial, atomic-facial and bisimplicial with respect to any linear functional $b$ whose 
minimum on the unit sphere occurs in the interior of the cone representing the identity of $W$.
Any linear extension of $W/\Theta$ is a shelling order on the facets of the associated sphere~$\Delta$.
\end{theorem}
The maximal cones of $\F_\Theta$ are the unions over $\Theta$-classes of the maximal cones of~$\F$.
If $\Theta$ and $\Phi$ are congruences such that $\Theta$ refines $\Phi$ then the lattice homomorphism associated to $\Phi$ 
factors through the homomorphism associated to $\Theta$.
In this case $\F_\Phi$ is refined by $\F_\Theta$.
The main shortcoming of Theorem~\ref{main} is that it gives no means of knowing when $\F_\Theta$ is the normal fan of a polytope.
It would be helpful to have a criterion for determining which of these fans are normal fans, particularly if the criterion were decisive
for the examples given later in the introduction.
The fan $\F_\Theta$ is not necessarily simplicial, but we give necessary and sufficient conditions on~$\Theta$ for $\F_\Theta$ to be 
simplicial (Proposition~\ref{simp cond}).

The Malvenuto-Reutenauer algebra is $\K[S_\infty]:=\bigoplus_{n\ge 0}\K[S_n]$ for a field $\K$, with a product which takes permutations 
$u\in S_p$ and $v\in S_q$ to the sum of all shuffles of $u$ and $v$ in $S_{p+q}$.
Loday and Ronco~\cite{LRorder} pointed out that this product can be expressed as the sum of the elements in a certain interval in weak order.
Furthermore they showed that the products on the algebra of planar binary trees and the algebra of non-commutative symmetric functions can be expressed 
as sums over intervals in Tamari lattices and Boolean algebras respectively.
The maps in~\cite{LRorder} relating these partial orders are $\eta$ and the descent map.

These facts are explained and generalized using lattice congruences.
A family of lattice congruences $\Theta_n$ on the weak order on the symmetric groups $S_n$ is called {\em translational} and/or 
{\em insertional} under certain conditions defined in Sections~\ref{trans sec} and~\ref{ins}.
Given any family $\set{\Theta_n}_{n\ge 0}$ of congruences, let $\set{Z^\Theta_n}_{n\ge 0}$ be the family of lattice quotients 
$S_n/\Theta_n$, and define a graded vector space $\K[Z^\Theta_\infty]:=\bigoplus_{n\ge 0}\K[Z^\Theta_n]$.
Define a map $c:\K[Z^\Theta_\infty]\to \K[S_\infty]$ by sending each element $x\in Z^\Theta_n$ to the sum of the elements of the 
corresponding congruence class in $S_n$.
We define a product on $\K[Z^\Theta_\infty]$ using each partial order $Z^\Theta_n$ in a manner analogous to Loday and Ronco's order-theoretic 
characterization of the shuffle product.
We also define a coproduct on $\K[Z^\Theta_\infty]$.
The main algebraic results of this paper are the following theorems.

\begin{theorem}
\label{subalgebra}
If $\set{\Theta_n}_{n\ge 0}$ is a translational family then the map $c$ embeds $\K[Z^\Theta_\infty]$ as a subalgebra of $\K[S_\infty]$.
\end{theorem}

\begin{theorem}
\label{subcoalgebra}
If $\set{\Theta_n}_{n\ge 0}$ is an insertional family then the map $c$ embeds $\K[Z^\Theta_\infty]$ as a subcoalgebra of $\K[S_\infty]$.
\end{theorem}
A translational and insertional family of congruences is called an {\em $\scrH$-family}, where the ``$\scrH$'' indicates ``Hopf,'' in accordance 
with the following immediate corollary of Theorems~\ref{subalgebra} and~\ref{subcoalgebra}.
\begin{corollary}
\label{subHopf}
If $\set{\Theta_n}_{n\ge 0}$ is an $\scrH$-family then the map $c$ embeds $\K[Z^\Theta_\infty]$ as a sub Hopf algebra of $\K[S_\infty]$.
\end{corollary}
The antipode of $\K[Z^\Theta_\infty]$ is easily written in terms of the antipode in $\K[S_\infty]$ (see Remark~\ref{anti remark}).
Given two $\scrH$-families $\set{\Theta_n}$ and $\set{\Phi_n}$ such that $\Theta_n$ refines $\Phi_n$ for each $n$, 
$\K[Z^\Phi_\infty]$ is a sub Hopf algebra of $\K[Z^\Theta_\infty]$.

The Tamari lattice is the subposet (in fact sublattice) of $S_n$ consisting of 312-avoiding permutations 
or alternately 231-avoiding permutations~\cite{Nonpure II}.
The left descent map on $S_n$ is a projection down to permutations avoiding both 231 and 312, and it follows that the descent map 
factors through~$\eta$.
Applying Theorem~\ref{main} recovers the refinement relationships on the associated fans.
The congruences associated to these lattices form $\scrH$-families, so Theorem~\ref{subHopf} can be applied.

The geometric results of this paper also apply to a broad generalization of the Tamari lattices, defined in~\cite{Cambrian}.
For any finite Coxeter group $W$ a family of {\em Cambrian congruences} is defined on the weak order on $W$.
The quotient of $W$ by a Cambrian congruence is called a {\em Cambrian lattice}.
The fans associated, via Theorem~\ref{main}, to the Cambrian congruences are conjectured to be combinatorially equivalent to the normal 
fans of the generalized associahedra~\cite{ga} and this conjecture is proven in types A and B\@.

For a general $\scrH$-family, a basis for each $\K[Z^\Theta_n]$ is characterized (Theorem~\ref{c-a avoid}) by a variation on pattern 
avoidance.\footnote{The pattern-avoidance description indicates that $\K[Z^\Theta_\infty]$ can also be obtained via an elegant general 
construction, due to Duchamp, Hivert, Novelli and Thibon, of sub Hopf algebras of $\K[S_\infty]$.
See Remark~\ref{monoid}.}
In the present paper, we exhibit several additional examples.
One of these examples builds $\K[S_\infty]$ as the limit of a sequence of smaller Hopf algebras $\K[S_{\infty,k}]$ where the first Hopf 
algebra in the sequence is a graded Hopf algebra with one-dimensional graded pieces (the binomial Hopf algebra $B_1$ 
of~\cite[Section V.2]{Joni-Rota}) and the second is the Hopf algebra of non-commutative symmetric functions.
Another example builds the Hopf algebra of planar binary trees from a similar sequence.

A third example concerns the {\em twisted Baxter permutations}, a set of permutations defined similarly to, and apparently equinumerous with, 
the Baxter permutations of~\cite{CGHK}.
By Theorem~\ref{c-a avoid}, the subposet of weak order on $S_n$
consisting of the twisted Baxter permutations is in fact the quotient of the weak order by a certain lattice congruence.
This congruence is identified as the meet of two congruences, one of which defines the Tamari lattice as the 
231-avoiding permutations, while the other defines the Tamari lattice as the 312-avoiding permutations.
The family of congruences defining the twisted Baxter permutations is an $\scrH$-family, so Theorem~\ref{subHopf} shows that there is a 
sub Hopf algebra of $\K[S_\infty]$ such that a basis for the $n$th graded piece is indexed by the twisted Baxter permutations in $S_n$.
Theorem~\ref{main}, besides proving several nice properties of the subposet of $S_n$ consisting of twisted Baxter permutations, also 
constructs a (non-simplicial) complete fan in $\reals^{n-1}$ whose maximal cones are indexed by the twisted Baxter permutations.
It would be interesting to know if this fan is the normal fan of some polytope.

This paper is the second in a series of papers beginning with~\cite{congruence} and continuing in~\cite{Cambrian}.
Each paper relies on the results of the preceding papers and cites later papers only for motivation or in the context of examples.

The organization of the remainder of this paper is as follows:
In Section~\ref{lat cong}, we provide background information on lattice congruences.
Section~\ref{fan posets} defines fan posets and exhibits their basic properties.
Section~\ref{po} defines the poset of regions of a central hyperplane arrangement, and quotes results which show that this poset is a fan poset 
with particularly nice properties.
In Section~\ref{quo} we prove a generalization of Theorem~\ref{main}.
Section~\ref{symmetric} provides background on the weak order on $S_n$ which is necessary for the exposition and proof of 
Theorem~\ref{subalgebra} in Section~\ref{trans sec} and Theorem~\ref{subcoalgebra} in Section~\ref{ins}.
Section~\ref{H} presents the characterization of $\scrH$-families by pattern avoidance, remarks on computing the product, coproduct 
and antipode in $\K[Z^\Theta_\infty]$.
The paper concludes with examples in Section~\ref{ex}.

\section{Lattice congruences}
\label{lat cong}
In this section we give background information on lattice congruences.
A more detailed exposition of lattice congruences can be found for example in~\cite{Gratzer}.
The poset notation used here is standard, and we assume basic poset and lattice terminology as for example in~\cite{EC1}.
If $x<y$ in $P$ and there is no $z\in P$ with $x<z<y$, say~$y$ {\em covers}~$x$ and write $x\covered y$.
If $P$ is a poset with a unique minimal element and a unique maximal element (for example if $P$ is a finite lattice), then the minimal element
is denoted $\0$ and the maximal element is $\1$.
The elements covering $\0$ are called the {\em atoms} of $P$, and the elements covered by $\1$ are {\em coatoms}.

Let $P$ be a finite poset with an equivalence relation $\Theta$ defined on the elements of $P$.
Given $a\in P$, let $[a]_\Theta$ denote the equivalence class of $a$.
The equivalence relation is an {\em order congruence} if:
\begin{enumerate}
\item[(i) ] Every equivalence class is an interval.
\item[(ii) ] The projection $\pidown:P\rightarrow P$, mapping each element $a$ of $P$ to the minimal element in $[a]_\Theta$, is order-preserving.
\item[(iii) ] The projection $\piup:P\rightarrow P$, mapping each element $a$ of $P$ to the maximal element in $[a]_\Theta$, is order-preserving.
\end{enumerate}

Define a partial order on the congruence classes by $[a]_\Theta\le[b]_\Theta$ if and only if there exists $x\in[a]_\Theta$ and $y\in[b]_\Theta$ 
such that $x\le_Py$.
The set of equivalence classes under this partial order is $P/\Theta$, the {\em quotient} of $P$ with respect to $\Theta$.
The quotient $P/\Theta$ is isomorphic to the induced subposet $\pidown(P)$.
The map $\piup$ maps $\pidown(P)$ isomorphically onto $\piup(P)$.
The inverse is $\pidown$.
For more information on order congruences and quotients, see~\cite{Cha-Sn,Order}.

The {\em join} $\join X$ of a subset $X\subseteq P$ is the unique minimal element, if it exists, of the set 
$\set{y\in P:y\ge x\mbox{ for all } x\in X}$.
Dually the {\em meet} $\meet X$ is the unique maximal lower bound of $X$.
A finite poset $L$ is called a {\em lattice} if every subset of $L$ has both a meet and a join.
An element $\gamma$ in a finite lattice $L$ is join-irreducible if and only if it covers exactly one element, which we denote $\gamma_*$.
The subposet of $L$ consisting of join-irreducible elements is denoted $\Irr(L)$.
A {\em lattice congruence} is an equivalence relation on a lattice which respects joins and meets.
Specifically, if $a_1\equiv a_2$ and $b_1\equiv b_2$ then $a_1\join b_1\equiv a_2\join b_2$ and similarly for meets.
When $L$ is a finite lattice, order congruences are exactly lattice congruences, and the quotient construction described above corresponds
to the algebraic notion of the quotient of a lattice with respect to a congruence.

The following simple properties\footnote{Quite likely these are known but they have not, to the author's knowledge, appeared in print.} of lattice congruences do not hold in the generality of poset congruences.
Let $\Theta$ be a congruence on a lattice $L$.
For $x\in L$, let $[x]_\Theta$ denote the congruence class of~$x$ mod $\Theta$.
\begin{lemma}
\label{int}
If $[x,y]$ is an interval in $L$, then $\set{[z]_\Theta:z\in[x,y]}$ is the interval $[[x]_\Theta,[y]_\Theta]$ in $L/\Theta$, and this interval is 
isomorphic to $[x,y]/\Theta$, where $\Theta$ also denotes the restriction of $\Theta$ to $[x,y]$.
\end{lemma}
\begin{proof}
If $z\in[x,y]$ then by definition $[z]_\Theta\in[[x]_\Theta,[y]_\Theta]$.
If $[z]_\Theta\in[[x]_\Theta,[y]_\Theta]$ then in particular $y\ge\pidown z$ and $x\le\piup z$.
So~$y$ and $\piup z$ are both upper bounds on~$x$ and $\pidown z$ and thus $x\join\pidown z$ is below both~$y$ and $\piup z$.
Thus we have $x\join\pidown z\in[x,y]$ and $x\join\pidown z\in[\pidown z,\piup z]=[z]_\Theta$, so 
$[z]_\Theta=[x\join\pidown z]_\Theta\in\set{[w]_\Theta:w\in[x,y]}$.

Since the interval $[x,y]$ is in particular a sublattice of $L$, the restriction of $\Theta$ to $[x,y]$ is a lattice congruence and the join and meet
operations in $[x,y]$ are inherited from $L$.
Therefore the join and meet operations in $L$ on congruence classes intersecting $[x,y]$ are the same as the join and meet of the restrictions 
of those congruence classes to $[x,y]$.
Thus $[[x]_\Theta,[y]_\Theta]$ and $[x,y]/\Theta$ are isomorphic as lattices.
\end{proof}

\begin{prop}
\label{quotient covers}
Let $L$ be a finite lattice, $\Theta$ a congruence on $L$ and $x\in L$.
Then the map $y\mapsto[y]_\Theta$ restricts to a one-to-one correspondence between elements of $L$ covered by $\pidown x$ and elements of 
$L/\Theta$ covered by $[x]_\Theta$.
\end{prop}
\begin{proof}
First, we show that the restriction of the map $y\mapsto[y]_\Theta$ to elements covered by $\pidown x$ is one-to-one.
Suppose that~$y$ and $y'$ are both covered by $\pidown x$ and $y\equiv y'$.
If $y\neq y'$ then $\pidown x$ is a minimal upper bound for~$y$ and $y'$, so it is in fact their join, and in particular $y\equiv \pidown x$.
This contradicts the fact that $\pidown x$ is the minimal element of its congruence class, thus proving that $y=y'$.
We now show that $[y]_\Theta\covered[x]_\Theta$ if and only if there is some $y'\equiv y$ such that $y'\covered\pidown x$ in $L$.

Suppose that $[y]_\Theta\covered[x]_\Theta$.
Thus in particular $\pidown y<\pidown x$, so let $y'$ be any element of $L$ such that $\pidown y\le y'\covered \pidown x$.
If $y'\not\equiv y$ then $[y]_\Theta<[y']_\Theta<[x]_\Theta$, which is a contradiction.
Thus the element $y'$ covered by $\pidown x$ has $[y']_\Theta=[y]_\Theta$.

Suppose that $y\covered\pidown x$ in $L$.
We want to show that $[y]_\Theta\covered[x]_\Theta$ in $L/\Theta$.
Since $\pidown x$ is minimal in $[\pidown x]_\Theta$ we have $y\not\equiv\pidown x$, so $[y]_\Theta<[x]_\Theta$.
Suppose that $[y]_\Theta\covered[z]_\Theta\le [x]_\Theta$ for some $z$ and let $y'$ be the unique element of
$[y]_\Theta$ covered by $\pidown z$ whose existence was proved in the previous paragraphs.
If $\pidown z\le y$, then $\pidown z\le \pidown y$, thus contradicting our supposition.
Since $\pidown x$ is an upper bound for~$y$ and $\pidown z$, we have $y\join\pidown z\le\pidown x$, and since $y\not\ge\pidown z$ and 
$\pidown x\covers y$, we have $y\join\pidown z=\pidown x$.
Now, since $y\equiv y'$, we have $y\join\pidown z\equiv y'\join\pidown z$, or in other words $\pidown x\equiv\pidown z$, so that in particular
$[z]_\Theta=[x]_\Theta$.
\end{proof}

Congruences on $L$ are, in particular, partitions of the elements of $L$, and $\Con(L)$ is the set of congruences of $L$
partially ordered by refinement.
The partial order $\Con(L)$ is a distributive lattice~\cite{Fun-Nak}, and thus is uniquely determined by the subposet $\Irr(\Con(L))$.
The meet in $\Con(L)$ is intersection of the congruences as relations.
If $\Theta_1$ and $\Theta_2$ are congruences on $L$, with associated downward projections $(\pidown)_1$ and $(\pidown)_2$, let 
$\Theta_1\join\Theta_2$ have associated downward projection $\pidown$.
It follows immediately from~\cite[Theorem I.3.9]{Gratzer} that $x\in L$ has $\pidown x=x$ if and only if both $(\pidown)_1 x=x$ and 
$(\pidown)_2 x =x$.
Thus the quotient of $L$ mod $\Theta_1\join\Theta_2$ is isomorphic to the induced subposet $((\pidown)_1L)\cap((\pidown)_2L)$ of $L$.

If $x\covered y$ and $x\equiv y\mod\Theta$, we say $\Theta$ {\em contracts} the edge $x\covered y$.
For an element~$y$, if there exists an edge $x\covered y$ contracted by $\Theta$, we say $\Theta$ contracts~$y$.
Thus $\Theta$ contracts a join-irreducible $\gamma$ if and only if $\gamma\equiv \gamma_*$.
A lattice congruence is determined by the set of join-irreducibles it contracts (see for example \cite[Section II.3]{Free Lattices}).
Given a covering pair $x\covered y$ in $L$, let $\Cg(x,y)$ be the smallest lattice congruence contracting that edge.
Then $\Cg(x,y)$ is a join-irreducible congruence.
Given a join-irreducible $\gamma$ of $L$, write $\Cg(\gamma)$ for $\Cg(\gamma_*,\gamma)$.
The map $\Cg:\Irr(L)\to\Irr(\Con(L))$ is onto, but need not be one-to-one.
A lattice $L$ is {\em congruence uniform} if $\Cg$ is a bijection and if a dual statement about meet-irreducibles holds as 
well~\cite{cong norm}.
When $L$ is a congruence uniform lattice, $\Irr(\Con(L))$ can be thought of as a partial order on the join-irreducibles.
If $\Theta$ is a congruence on $L$, then $\Irr(\Con(L/\Theta))$ is the order filter in $\Irr(\Con(L))$ consisting of join-irreducibles of $L$ not 
contracted by $\Theta$.

Given a congruence $\Theta_1$ on a lattice $L_1$ and a congruence $\Theta_2$ on a lattice $L_2$, define an equivalence 
$\Theta_1\times\Theta_2$ on $L_1\times L_2$ by setting $(x_1,x_2)\equiv(y_1,y_2)$ mod $\Theta_1\times\Theta_2$ if and only if 
$x_1\equiv y_1$ mod $\Theta_1$ and $x_2\equiv y_2$ mod $\Theta_2$.
It is an easy exercise to show that $\Theta_1\times\Theta_2$ is a congruence, and furthermore that any congruence on $L_1\times L_2$ 
has the form $\Theta_1\times\Theta_2$ for some congruence $\Theta_1$ on $L_1$ and some congruence $\Theta_2$ on $L_2$.
The join-irreducibles of $L_1\times L_2$ are exactly the pairs $(\gamma_1,\0)$ where $\gamma_1$ is a join-irreducible of $L_1$, and 
the pairs $(\0,\gamma_2)$ where $\gamma_2$ is a join-irreducible of $L_2$.

Given lattices $L_1$ and $L_2$ a {\em homomorphism} from $L_1$ to $L_2$ is a map $\eta:L_1\to L_2$ such that for all~$x$ and~$y$ in 
$L_1$ we have $\eta(x\join y)=\eta(x)\join\eta(y)$ and similarly for meets.
Given a lattice homomorphism $\eta$, the equivalence relation whose classes are the fibers of $\eta$ is a congruence,
and conversely, given a congruence $\Theta$ on $L$, the map from an element to its equivalence class is a homomorphism $L\to(L/\Theta)$.
Alternately, the map $\pidown$ is a homomorphism from $L$ to $\pidown L\cong L/\Theta$.
If $\eta_1:L\to L_1$ and $\eta_2:L\to L_2$ are lattice homomorphisms, we say $\eta_2$ {\em factors through} $\eta_1$ if there is a lattice
homomorphism $\eta:L_1\to L_2$ such that $\eta_2=\eta\circ\eta_1$.
If $\Theta_1$ and $\Theta_2$ are the lattice congruences associated to $\eta_1$ and $\eta_2$ and $\Theta_1\le\Theta_2$ in $\Con(L)$ 
then $\eta_2$ factors through $\eta_1$.

Given a partially ordered set $P$, topological statements about $P$ refer to its {\em order complex}, the abstract simplicial complex whose 
faces are the {\em chains} (totally ordered subposets) of $P$.
The {\em proper part} of a finite lattice $L$ is $L-\set{\0,\1}$.
The following is a special case of the Crosscut Theorem (see the explanation surrounding (10.8) of~\cite{Top}).
\begin{theorem}
\label{crosscut}
If $L$ is a finite lattice with atoms $A$, then the proper part of $L$ is homotopy equivalent to the abstract simplicial complex consisting of 
subsets of $A$ whose join is not $\1$.
\end{theorem}
For convenience here, we call this abstract simplicial complex the {\em crosscut complex} of $L$, although the usual definition of a crosscut 
complex is much more general.

\begin{corollary}
\label{cross cor}
If $L$ is a lattice and $\Theta$ is a congruence on $L$ such that no atom of $L$ is congruent to $\0$ and no coatom is congruent to $\1$, then 
the proper part of $L$ is homotopy equivalent to the proper part of $L/\Theta$.
\end{corollary}
\begin{proof}
Since no atom of $L$ is congruent to $\0$, there is a one-to-one correspondence between atoms of $L$ and atoms of $L/\Theta$.
We use $A$ to denote both sets of atoms.
Since no coatom is congruent to $\1$, the top element of $L/\Theta$ is the equivalence class $\set{\1}$.
For $S\subseteq A$ the join of $A$ in $L/\Theta$ is the equivalence class of the join of $A$ in $L$, so $S$ joins to $\1$ in $L$ if and only if 
it joins to $\set{\1}$ in $L/\Theta$.
Thus the crosscut complex of $L/\Theta$ is isomorphic to the crosscut complex of $L$.
\end{proof}

\section{Fan Posets}
\label{fan posets}
In this section we define fan posets and prove some of their basic properties.
We assume the definitions of polytopes, cones, simplicial complexes, regular CW complexes, combinatorial isomorphism and homotopy equivalence.
For more information on regular CW complexes, particularly as they relate to combinatorics, 
see~\cite{Top} and Section 4.7 of~\cite{Oriented Matroids}.
We call the closed cells of a CW complex {\em faces}. 
The {\em 1-skeleton} of a CW complex $\Gamma$ is the subcomplex consisting of the 0- and 1-dimensional faces of $\Gamma$.
Given a regular CW complex $\Gamma$ with face poset $P$, the poset $P-\set{\0}$ is topologically equivalent to $\Gamma$, because the order
complex of $P-\set{\0}$ is combinatorially isomorphic to the {\em barycentric subdivision} of $\Gamma$.
The following theorem is due to Bj\"{o}rner~\cite{Bj}.
\begin{theorem}
\label{CW thm}
A non-trivial poset $P$ with a unique minimal element $\0$ is the face poset of a regular CW complex if and only if every  interval $(\0,x)$ 
is a sphere.
\end{theorem}

A {\em fan} in $\reals^d$ is a family $\F$ of nonempty closed polyhedral cones with the following properties:
\begin{enumerate}
\item[(i) ] Every nonempty face of a cone in $\F$ is also a cone in $\F$.
\item[(ii) ] The intersection of two cones in $\F$ is a face of both.
\end{enumerate}
A {\em complete fan} has the additional property that $\cup\F=\reals^d$.
Since $\F$ is closed under intersections and has a minimal element $\cap\F$, if one partially orders $\F$ by inclusion and adjoins a maximal 
element $\1$, one obtains a lattice, called the {\em face lattice} of $\F$.
The intersection $\cap\F$ of all cones in $\F$ is a subspace, because otherwise it must have at least one proper nonempty face.
If for every maximal cone $C$ of $\F$, the normals to the facets of $C$ are linearly independent, then $\F$ is a {\em simplicial fan}.
We need the following lemma.
\begin{lemma}
\label{fan lemma}
Let $\C$ be a finite set of $d$-dimensional closed cones in $\reals^d$ with non-intersecting interiors such that $\cup\C=\reals^d$, with the 
property that the intersection of any two cones in $\C$ is a face of each.
Then the collection $\F$ of cones, consisting of arbitrary intersections of cones in $\C$, is a fan.
\end{lemma}
\begin{proof}
We first show that $\F$ is exactly the set of faces of cones in $\C$.
Let $F$ be the intersection of some subset $S\subseteq\C$.
We prove by induction on $|S|$ that $F$ is a face of some $C\in S$.
The base of the induction is the case $|S|=1$, or in other words $F\in\C$, which is trivial.
Let $S'=S-\set{D}$ for some $D\in S$.
Then by induction, $\cap S'$ is a face of some $C\in S'$.
Since $C\cap D$ is a face of $C$ as well, the intersection $\cap S$ is the intersection of two faces of $C$, and thus is a face of $C$.

Conversely, let $G$ be a face of some cone $C$ in $\C$.
Then $G$ can be written as the intersection of some set $M$ of facets of $C$.
Since $\cup\C=\reals^d$, and the members of $\C$ intersect in faces, each facet $F$ of $C$ is the intersection of $C$ with some $C_F\in\C$.
Then $G$ is the intersection of $C$ with all of the $C_F$ for $F$ in $M$.

We have thus established condition (i) in the definition of a fan.
Furthermore, since each cone in $\F$ is the intersection of some set of cones in $\C$, the intersection of two cones $F$ and $G$ in $\F$ is also 
the intersection of some set of cones in $\C$, and thus $F\cap G$ is the face of some cone in $\C$.
Since $F$ and $G$ are each faces of cones in $\C$, the intersection $F\cap G$ is a face of each.
\end{proof}

Given a cone $C$ of $\F$, we define the restriction $\F|_C$ of $\F$ to $C$ as follows.
Let $U$ be an open $d$-ball centered at a point $p$ in the relative interior of $C$, such that $U$ does not intersect any cone not containing $C$.
Then the intersection of $\F$ with $U$ gives a cellular decomposition of $U$.
We center a $d$-dimensional vector space at $p$ and extend this cellular decomposition of $U$ linearly to a cellular decomposition of the 
vector space.
The resulting decomposition is the fan $\F|_C$.

A complete fan $\F$ is {\em essential} if $\cap\F$ is the origin.
If $\F$ is essential then the intersection of $\F$ with the unit sphere defines a cellular decomposition of the sphere.
Given a non-essential fan $\F$, a combinatorially isomorphic essential fan $\F/(\cap\F)$ is obtained by intersecting $\F$ with the 
orthogonal complement $(\cap\F)^\perp$ of $\cap\F$.
We define the {\em associated sphere} $\Delta$ of $\F$ to be the CW sphere whose cellular structure is the decomposition of the unit sphere in 
$(\cap\F)^\perp$ induced by $\F/(\cap\F)$.
The upper interval $[C,\1]$ in the face lattice of $\F$ is isomorphic to the face lattice of $\F|_C$.
In particular, $[C,\1]$ is spherical, and by Theorem~\ref{CW thm} there is a regular CW sphere $\Gamma$ whose face poset, upon adjoining a
maximal element $\1$, is dual to the face lattice of $\F$.
We call $\Gamma$ the {\em dual sphere} to $\F$.

If $\Pi$ is a polytope and $F$ is a face of $\Pi$, the {\em normal cone} to $F$ is the set of linear functionals which are 
maximized at every point on $F$.
The {\em normal fan} of $\Pi$ is the collection of normal cones to the faces of $\Pi$.
A fan is called {\em polytopal} if it is the normal fan of some polytope.

A {\em fan poset} is a pair $(\F,P)$ where $\F$ is a complete fan in $\reals^d$ and $P$ is a finite poset whose elements are the 
maximal cones of $\F$, subject to the following conditions:
\begin{enumerate}
\item[(i) ]For every interval $I$ of $P$, the union of the maximal cones in $I$ is a polyhedral cone.
\item[(ii) ]For every cone $C$ of $\F$, the set of maximal cones containing $C$ is an interval in $P$.
\end{enumerate}
The intervals arising as in (ii) are called {\em facial} intervals.

Say $(\F,P)$ is {\em homotopy facial} if the homotopy types of intervals are described as follows:
if $[x,y]$ is a facial interval associated to a cone of dimension $k$, then the open interval $(x,y)$ is homotopy equivalent to a 
$(d-2-k)$-sphere.
If $[x,y]$ is not a facial interval, then $(x,y)$ is contractible.
By convention the complex containing only the empty set is a $(-1)$-dimensional sphere, and the empty complex is a $(-2)$-dimensional 
sphere.
If $(\F,P)$ is homotopy facial then in particular the face lattice of $\F$ can be determined from the abstract partial order $P$.
It is dual to the set of non-contractible intervals, partially ordered by containment.
The non-contractible intervals in $P$ are exactly the intervals with non-zero M\"{o}bius functions.

An interval $I$ in a poset is called {\em atomic} if the maximal element of $I$ is the join of the set of atoms of $I$.
Call $(\F,P)$ {\em atomic-facial} if the facial intervals are exactly the atomic intervals.
If $(\F,P)$ is atomic facial, then the face lattice of $\F$ is dual to the set of atomic intervals, partially ordered by containment.

Let $(\F,P)$ be a fan poset and let $b$ be a linear functional on $\reals^d$.
For any covering relation $C_1\covered C_2$ in $P$, let $\nu$ be the unit normal vector to the hyperplane separating $C_1$ from $C_2$, 
oriented to point from $C_1$ to $C_2$. 
Say that $(\F,P)$ is {\em induced by} $b$ if for any such $C_1\covered C_2$ and $\nu$ we have $b(\nu)>0$.
For any maximal cone $C$ of $\F$, let $N_+$ be the set of outward-facing unit normals $\nu$ of $C$ such that $b(\nu)>0$ and let $N_-$ 
be the set of outward-facing unit normals for which $b(\nu)<0$.
Say a maximal cone $C$ is {\em bisimplicial} with respect to $b$ if both $N_+$ and $N_-$ are linearly independent sets.
Say $(\F,P)$ is {\em bisimplicial} with respect to $b$ if it induced by $b$ and if each maximal cone of $\F$ is bisimplicial with respect to $b$.

A facial interval of a fan poset is itself a fan poset.
That is, if $C$ is a cone of $\F$ and $I$ is the corresponding interval of $P$, then $(\F|_C,I)$ is a fan poset.
If $(\F,P)$ is polytopal, homotopy facial, atomic-facial, induced, bisimplicial and/or simplicial, then $(\F|_C,I)$ enjoys those properties as well.

The fan poset $(\F,P)$ is defined to be {\em polytopal} and/or {\em simplicial} if $\F$ is.
If $(\F,P)$ is polytopal, then the polytope is combinatorially isomorphic to $\Gamma$, so we refer to the polytope as $\Gamma$.
If $(\F,P)$ is a polytopal fan poset induced by a linear functional~$b$, then $P$ is the partial order induced by $b$ on the vertices of $\Gamma$.

Suppose $P$ is a partial order on the vertices of some CW sphere $\Gamma$.
Say $P$ {\em orients the 1-skeleton} of $\Gamma$ if the 1-skeleton of $\Gamma$ is isomorphic as a 
graph to the Hasse diagram of $P$ via the identification of elements of $P$ with vertices of $\Gamma$.
\begin{prop}
\label{1skel}
If $(\F,P)$ is a fan poset then $P$ orients the 1-skeleton of the dual sphere $\Gamma$.
\end{prop}
\begin{proof}
Edges in $\Gamma$ correspond to pairs of maximal cones of $\F$ intersecting in dimension $d-1$.
Condition (ii) in the definition of a fan means in particular that the 1-skeleton of $\Gamma$ has no multiple edges.
Thus showing graph isomorphism is equivalent to showing that two maximal cones form a cover in $P$ if and only if the maximal cones 
intersect in dimension $d-1$.

Suppose $C_1$ and $C_2$ are maximal cones of $\F$ such that $F:=C_1\cap C_2$ has dimension $d-1$.
Then $\set{C_1,C_2}$ is the complete set of maximal cones containing $F$.
By the definition of a fan poset, $\set{C_1,C_2}$ is an interval in $P$, necessarily a cover relation.

Suppose $C_1\covered C_2$ in $P$, so that in particular $\set{C_1,C_2}$ is an interval in $P$.
By the definition of fan poset, $C_1\cup C_2$ is a polyhedral cone, so in particular $C_1$ and $C_2$ must intersect in dimension $d-1$.
\end{proof}

Let $\Delta$ be a CW complex all of whose facets have dimension $d$.
Let $\delta F$ denote the boundary of a face (closed cell) $F$ of $\Delta$.
A linear order $F_1,F_2,\ldots,F_t$ on the facets of $\Delta$ is a {\em shelling} of $\Delta$ if $d=0$ or if $d\ge 1$ and the following conditions 
hold:
\begin{enumerate}
\item[(i) ] $\delta F_1$ has a shelling,
\item[(ii) ] For $2\le j\le t$, the intersection $F_j\cap(\cup_{i=1}^{j-1}F_j)$ is a pure CW complex of dimension $d-1$, and
\item[(iii) ] For $2\le j\le t$, the boundary $\delta F_j$ has a shelling in which the $(d-1)$-dimensional faces of 
$F_j\cap(\cup_{i=1}^{j-1} F_j)$ appear first.
\end{enumerate}

The {\em boundary complex} of a convex polytope is the boundary of the polytope with a cellular decomposition consisting of the relative interiors
of the faces of the polytope.
Bruggesser and Mani~\cite{BM} defined, for any linear functional $b$ not parallel to any facet hyperplane of the polytope, a shelling of 
the boundary complex.
Their shelling has the following property: 
Every facet whose outward-facing normal $\nu$ has $b(\nu)<0$ precedes every facet whose outward-facing normal $\nu$ has $b(\nu)>0$.
\begin{prop}
\label{induced shelling}
If $(\F,P)$ is a fan poset induced by a linear functional then any linear extension of $P$ is a shelling order on the associated sphere $\Delta$.
\end{prop}
\begin{proof}
Let $\Delta$ be the $d'$-dimensional sphere associated to $\F$ and let $F_1,F_2,\ldots,F_t$ be a linear order on the facets of $\Delta$ 
induced by some linear extension of $P$.
For each~$i$ let $C_i$ be the maximal cone of $\F$ containing $F_i$.
We first establish  condition (ii) in the definition of shelling, independent of the hypothesis that $(\F,P)$ is induced by a linear functional.
Suppose $1\le i<j\le t$.
The maximal cones of $\F$ containing $C_i\cap C_j$ form an interval in $P$.
Since $i<j$, in particular $C_j$ is not the bottom element of the interval, so we can find a maximal cone $C_k$ which is in the interval and is 
covered by $C_j$.
Since $C_k$ is covered by $C_j$, by the proof of Proposition~\ref{1skel} the intersection $C_j\cap C_k$ is a facet of $C_j$, and since 
$C_k$ is in the interval, we have $C_i\cap C_j\subseteq C_k$.
Finally, since $C_k$ is below $C_j$ in $P$, we have $k<j$.
Intersecting with the unit sphere in $(\cap\F)^\perp$, we have the following statement about $\Delta$:
For every $1\le i<j\le t$ there exists $1\le k<j$ such that $F_i\cap F_j\subseteq F_k$ and $F_j\cap F_k$ is a facet of $F_j$.
Thus every face of $F_j\cap(\cup_{i=1}^{j-1}\delta F_j)$ is contained in a $(d'-1)$-dimensional face of 
$\delta F_j\cap(\cup_{i=1}^{j-1}\delta F_j)$, implying condition (ii) in the definition of shelling.

We can assume $\F$ is essential because if not, we replace $\F$ by $\F/(\cap\F)$.
Since $F_1,\ldots,F_t$ is a linear extension of $P$, the set $C_j\cap(\cup_{i=1}^{j-1}C_j)$ is exactly the union of the facets of $C_j$ which 
separate $C_j$ from maximal cones covered by $C_j$ in $P$.
If $(\F,P)$ is induced by $b$ then this set of facets is exactly the set of facets whose outward-facing normals $\nu$ have $b(\nu)<0$.
Intersecting $C_j$ with an affine hyperplane $H$ parallel to $b$ so as to produce a convex polytope of dimension $d-1$, the 
Bruggesser-Mani shelling with respect to $b$ is a shelling of $C_j\cap H$ in which the $(d-2)$-dimensional faces of 
$C_j\cap H\cap(\cup_{i=1}^{j-1}C_j)$ appear first.
Since $C_j\cap H$ is combinatorially isomorphic to $F_j$, this satisfies (iii).
We can shell $F_1$ in a similar manner, using any linear functional not parallel to a facet of $F_1$. 
\end{proof}
Our proof is patterned after the proof, due to Bj\"{o}rner and Ziegler, of a similar statement~\cite[Proposition 4.3.2]{Oriented Matroids} due to 
Lawrence about the ``big'' face poset of an oriented matroid.

Any linear ordering of the facets of a simplex is a shelling order.
Thus if $\Delta$ is a pure simplicial complex, a total order on the facets of $\Delta$ is a shelling if and only if it satisfies condition (ii) in the 
definition of shelling given above.
In the proof of Proposition~\ref{induced shelling}, condition (ii) was established independent of the condition that $(\F,P)$ is induced.
Thus we have the following:
\begin{prop}
\label{shelling}
If $(\F,P)$ is a simplicial fan poset, then any linear extension of $P$ is a shelling order on $\F$.
\end{prop}

One application of a shelling order on a simplicial complex is in determining the face numbers of the simplicial sphere $\Delta$ associated to 
a simplicial fan $\F$.
The {\em $f$-vector} of a simplicial complex $\Lambda$ of dimension $d-1$ is $(f_{-1},f_0,f_1,\ldots,f_{d-1})$, where $f_i$ is the number of 
simplices of $\Lambda$ of dimension $i$ and the empty simplex is by convention $(-1)$-dimensional.
The {\em $h$-vector} of $\Lambda$ is $(h_0,h_1,\ldots,h_d)$, defined by the polynomial identity
\[\sum_{i=0}^df_{i-1}(x-1)^{d-i}=\sum_{i=0}^dh_ix^{d-i}.\]
For example $f_0=h_1+d$ and $h_d$ is $(-1)^{d-1}$ times the reduced Euler characteristic of $\Lambda$.
When $\Lambda$ is shellable, for each maximal simplex $F_j$ in the shelling order there is a unique minimal face $\R(F_j)$ of $F_j$ among faces of $F_j$ 
not contained in $\cup_{i=1}^{j-1}F_i$.
Furthermore $\sum_{i=0}^dh_ix^i=\sum_{i=0}^tx^{|\R(F_i)|}$, where $|\R(F_i)|$ is the number of vertices of $\R(F_i)$.
Equivalently, $|\R(F_i)|$ is the number of facets of $F_i$ contained in $\cup_{i=1}^{j-1}F_i$.
For $(\F,P)$ a simplicial fan poset, $C$ a maximal cone of $\F$ and any linear extension of $P$, the quantity $|\R(C)|$ is the number of 
elements covered by $C$ in $P$.
Thus for $i=0,1,\ldots,d$ the number of elements of $P$ covering exactly $i$ elements is $h_i$ in the $h$-vector of $\Delta$.
This is in keeping with the fact that $P$ is a {\em good orientation} of $\Gamma$ in the sense of Kalai~\cite{Kalai}.

The Dehn-Sommerville equations $h_i=h_{d-i}$ for $i=0,1,\ldots,d$ are satisfied by the boundary complexes of simplicial polytopes.
If $(\F,P)$ is a simplicial fan poset and if $P'$ is the dual partial order to $P$, then $(\F,P')$ is a simplicial fan poset with the same 
associated sphere.
An element covering $i$ elements in $P$ covers $d-i$ elements in $P'$.
Since the $h$-vector is a combinatorial invariant of $\Delta$ we have the following.
\begin{prop}
\label{DS}
If $(\F,P)$ is a simplicial fan poset then the associated simplicial sphere satisfies the Dehn-Sommerville equations.
\end{prop}

For a cone $C$ of $\F$, the {\em star} of $C$ is the fan whose maximal cones are the maximal cones of $\F$ which contain $C$.
The star of $C$ is {\em convex} if the union of the maximal cones of the star is a convex set.
A fan $\F$ is {\em locally convex} if the star of every cone of $\F$ is convex.
The following is immediate from the definition of a fan poset.

\begin{prop}
\label{local convex}
If $(\F,P)$ is a fan poset then $\F$ is locally convex.
\end{prop}

A simplicial complex $\Delta$ is {\em flag} if every minimal set of vertices not spanning a face of $\Delta$ has cardinality 2.
In~\cite{Le-Rei} it is shown that, given a simplicial locally convex fan $\F$, the simplicial sphere $\Delta$ is flag.
Thus Proposition~\ref{local convex} implies the following.
\begin{prop}
\label{flag}
If $(\F,P)$ is a simplicial fan poset then the corresponding simplicial sphere is flag.
\end{prop}

\section{Posets of Regions}
\label{po}
In this section we give background information on the poset of regions of a hyperplane arrangement, prove or quote basic results, and restate some known results in the language
of fan posets.
The poset of regions was defined by Edelman~\cite{Edelman} and further studied in~\cite{BEZ,Ed-Wa,hyperplane,hplanedim}.

A {\em hyperplane arrangement}~$\A$ is a finite collection of codimension~1 linear subspaces in $\reals^d$ called {\em hyperplanes}.
The complement of the union of the hyperplanes is disconnected, and the closures of its connected components are called {\em regions}.
In general, one might consider an arrangement of affine hyperplanes.
Hyperplane arrangements consisting entirely of linear subspaces are called {\em central}, and all hyperplane arrangements considered in this 
paper are central.
The {\em rank} of an arrangement is the dimension of the linear span of the normals to the hyperplanes.
A region $R$ of~$\A$ is called {\em simplicial} if the normals to the facets of $R$ are linearly independent.
A central hyperplane arrangement is called {\em simplicial} if every region is simplicial.

We fix once and for all a central hyperplane arrangement~$\A$ and a region $B$ of~$\A$.
A hyperplane $H$ is said to {\em separate} two distinct points $x_1$ and $x_2$ in $\reals^d$ if the line segment whose endpoints are
$x_1$ and $x_2$ intersects $H$ in exactly one point.
For regions $R_1$ and $R_2$, a hyperplane $H\in\A$ {\em separates} $R_1$ from $R_2$ if $H$ separates any (or equivalently {\bf every})
pair of points $(x_1,x_2)$ with $x_1$ in the interior of $R_1$ and $x_2$ in the interior of $R_2$.

For any region $R$, define the {\em separating set} $S(R)$ of $R$ to be the set of hyperplanes separating $R$ from $B$. 
The {\em poset of regions} $\Po(\A,B)$ is a partial order on the regions with $R_1\le R_2$ if and only if $S(R_1)\subseteq S(R_2)$.
The region $B$, called the {\em base region}, is the unique minimal element of $\Po(\A,B)$.
The map sending each region $R$ to its antipodal region $-R$ is an anti-automorphism and corresponds to complementation of separating sets.
In particular, $\Po(\A,B)$ has a unique maximal element $-B$.
Given a region $R$, call those facets of $R$ by which one moves up in $\Po(\A,B)$ {\em upper facets} of $R$, and call the other facets of $R$
{\em lower facets}.

Associated to~$\A$ there is a complete fan which we call $\F$, consisting of the regions of~$\A$ together with all of their faces.
Given a cone $C$ of $\F$, the set of regions containing $C$ is an interval in $\Po(\A,B)$, isomorphic to the poset of regions $\Po(\A',B')$, where
$\A'$ is the set of hyperplanes of~$\A$ containing $C$ and $B'$ is the region of $\A'$ containing $B$.
Given an interval $[R_1,R_2]$ in $\Po(\A,B)$, the union of the corresponding regions is the closure of the set of points separated from $B$ 
by every hyperplane in $S(R_1)$ and separated from $-B$ by every hyperplane in $\A~-~S(R_2)$.
This set is a polyhedral cone, and thus $(\F,\Po(\A,B))$ is a fan poset.
The definition of $\Po(\A,B)$ by separating sets can be rephrased as the statement that $\Po(\A,B)$ is the partial order induced on the 
maximal cones of $\F$ by any linear functional $b$ whose minimum on the unit sphere lies in the interior of $B$.
The fan $\F$ is the normal fan to a zonotope which is the Minkowski sum of the normal vectors to the hyperplanes.
The dimension of the zonotope is the rank of the arrangement.
In particular, $(\F,\Po(\A,B))$ is polytopal.
Edelman and Walker~\cite[Theorem 2.2]{Ed-Wa} determined the homotopy type of open intervals in $\Po(\A,B)$.
In the terminology of fan posets, their theorem is exactly the statement that $(\F,\Po(\A,B))$ is homotopy facial.

\begin{lemma}
\label{joins}
If $\R$ is the set of regions covering $B$ in $\Po(\A,B)$ then $\join\R=-B$ and any proper subset $S\subsetneq\R$ has an upper bound strictly
below $-B$.
\end{lemma}
\begin{proof}
For any $R\in\R$, there is some $H\in\A$ such that $S(R)=\set{H}$ and $S(-R)=\A-\set{H}$.
Any element covered by $-B$ is $-R$ for some $R\in\R$, and in particular, no element covered by $-B$ is above every element of $\R$, so 
$\join\R=-B$.
For any $U\subsetneq\R$, take $R\in\R-U$, and let $H$ have $S(R)=\set{H}$.
Then $S(-R)=\A-\set{H}$ in particular contains $S(R')$ for every $R'\in U$, so $-R$ is an upper bound for $U$.
\end{proof}
If $I$ is a facial interval of $(\F,\Po(\A,B))$, then since $I$ is isomorphic to some other poset of regions, by Lemma~\ref{joins} it is an 
atomic interval of $\Po(\A,B)$.
If~$\A$ is simplicial and $I$ is an atomic interval, let $R$ be the minimal element of $I$, and let $A$ be the set of atoms of $I$.
Then $C:=R\cap(\cap A)$ is a face of $R$ and thus a cone $\F$.
Let $\A'$ be the set of hyperplanes of~$\A$ containing $C$.
The join of $A$ is the region containing $C$ whose separating set is $S(R)\cup\A'$, and thus $I$ is the complete set of regions containing $C$.

We summarize these facts in the following theorem:
\begin{theorem}
\label{poset}
If~$\A$ is a central hyperplane arrangement, $B$ is a region of~$\A$ and $\F$ is the corresponding fan, then 
\begin{enumerate}
\item[(i) ]$(\F,\Po(\A,B))$ is a fan poset, 
\item[(ii) ]$(\F,\Po(\A,B))$ is polytopal, homotopy facial, and induced by any linear functional $b$ whose minimum on the unit sphere lies in 
the interior of $B$.
\item[(iii) ]Facial intervals are atomic, and if~$\A$ is simplicial then $(\F,\Po(\A,B))$ is atomic-facial.
\end{enumerate}
\end{theorem}

The following easy lemma will be useful in a later section.
\begin{lemma}
\label{face}
For any two regions $Q$ and $R$ of~$\A$, there is a sequence of regions $Q=R_0,\ldots R_t=R$ such that for every $i$ the intersection 
$R_i\cap R_{i-1}$ is $(d-1)$-dimensional and $Q\cap R\subseteq R_i$ for every $i$.
\end{lemma}
\begin{proof}
We may as well take $Q=B$.
Then because $\F$ is a fan poset, the set of regions containing $Q\cap R$ is an interval in $\Po(\A,Q)$, and the desired sequence is any 
unrefinable chain from $Q$ to $R$ in the interval.
\end{proof}

Bj\"{o}rner, Edelman and Ziegler~\cite{BEZ} showed that if~$\A$ is simplicial, then $\Po(\A,B)$ is a lattice for any choice of $B$.
In~\cite{congruence} it is shown that when~$\A$ is simplicial then $\Po(\A,B)$ admits special congruences called {\em parabolic congruences},
which we now define.
Let~$\A$ be simplicial and let $\B$ be the set of facet hyperplanes of $B$, and for each $H\in\B$ let $R(H)$ be the atom of $\Po(\A,B)$ 
separated from $B$ by $H$.
For any $K\subseteq\B$ the intersection of the hyperplanes of $\B-K$ is a subspace $L_K$.
Let $\A_K$ be the set of hyperplanes containing $L_K$ and let $B_K$ be the $\A_K$-region containing $B$.
The arrangement $\A_K$ is simplicial.
Let $\Theta_K$ be the equivalence relation on $\Po(\A,B)$ setting $R_1\equiv R_2$ if and only if $R_1$ and $R_2$ are contained in the same
$\A_K$-region.
In other words, the $\A_K$ regions are the unions over $\Theta_K$-classes of the~$\A$-regions.
The equivalence $\Theta_K$ is a lattice congruence~\cite[Proposition 6.3]{congruence}, and furthermore:
\begin{theorem}~\cite[Theorem 6.9]{congruence}
\label{parabolic} 
Let~$\A$ be simplicial and let $K\subseteq \B$.
Then $\Theta_K$ is the unique minimal lattice congruence with $B\equiv R(H)$ for every $H\in(\B-K)$.
\end{theorem}
When~$\A$ is a Coxeter arrangement, the homomorphism associated to $\Theta_K$ is projection to a parabolic subgroup.

In the next section, for any congruence $\Theta$ on a lattice $\Po(\A,B)$ we construct a fan $\F_\Theta$ whose maximal cones are the unions
over $\Theta$-classes of the maximal cones of $\F$.
Suppose~$\A$ is simplicial, $K\subseteq\B$ and $\Theta$ is any congruence contracting atoms $R(H)$ for $H\in K$.
By Theorem~\ref{parabolic}, $\Theta$ is refined by $\Theta_K$, so that $\Theta$ can be thought of as a congruence on $\Po(\A_K,B_K)$.
Thus we can first pass to the fan associated to $\Po(\A_K,B_K)$ and form the fan $\F_\Theta$ by taking unions of $\A_K$-regions.
In particular, when~$\A$ is simplicial we can always reduce to the case where $\Theta$ contracts no atoms of $\Po(\A,B)$.
Furthermore, we have the following:
\begin{prop}
\label{may as well}
Let~$\A$ be simplicial and let $\Theta$ be a congruence on $\Po(\A,B)$.
If $[B]_\Theta\neq\set{B}$ then $(\cap\F)\subsetneq(\cap\F_\Theta)$.
\end{prop}
\begin{proof}
We have $\cap\F=\cap\A$.
If $[B]_\Theta\neq\set{B}$ then for some nonempty $K\subseteq\B$ we have $R(H)\equiv B$ for every $H\in K$.
Thus $(\cap\A)\subsetneq(\cap\A_K)\subseteq(\cap\F_\Theta)$.
\end{proof}
It follows easily from the definition of a lattice congruence that $[B]_\Theta=\set{B}$ if and only if $[-B]_\Theta=\set{-B}$.

\section{Congruences and fan lattices}
\label{quo}
This section is devoted to proving a generalization of Theorem~\ref{main} and other facts about the fans $\F_\Theta$.

\begin{theorem}
\label{main gen}
If~$\A$ is a central hyperplane arrangement and $B$ is a region of~$\A$ such that $\Po(\A,B)$ is a lattice, then for any lattice 
congruence $\Theta$ on $\Po(\A,B)$ there is a complete fan $\F_\Theta$, refined by $\F$, with the following properties:
\begin{enumerate}
\item[(i) ]$(\F_\Theta,\Po(\A,B)/\Theta)$ is a fan lattice.
\item[(ii) ]$(\F_\Theta,\Po(\A,B)/\Theta)$ is induced by any linear functional whose minimum on the unit sphere lies in the interior of $B$.
\item[(iii) ]Any linear extension of $\Po(\A,B)/\Theta$ is a shelling of $\F_\Theta$.
\item[(iv) ]If~$\A$ is simplicial then $(\F_\Theta,\Po(\A,B)/\Theta)$ is homotopy facial, atomic-facial and bisimplicial with respect to any 
linear functional whose minimum on the unit sphere lies in the interior of $B$.
\end{enumerate}
\end{theorem}
If~$\A$ is a Coxeter arrangement then~$\A$ is simplicial and $\Po(\A,B)$ is a lattice isomorphic to the weak order on the associated 
Coxeter group.
Thus Theorem~\ref{main} is a special case of Theorem~\ref{main gen}.
As mentioned in the introduction, when $(\F_\Theta,\Po(\A,B)/\Theta)$ is homotopy facial, the M\"{o}bious function of a non-facial interval 
in $\Po(\A,B)/\Theta$ is zero, and a facial interval for a face $F$ has M\"{o}bius function $(-1)^{d-2-\dim F}$.

To construct the fan $\F_\Theta$ for $\Po(\A,B)/\Theta$, recall that each congruence class of~$\Theta$ is an interval in $\Po(\A,B)$,
 so the union of the corresponding maximal cones of~$\F$ is a convex cone.
Let $\C$ be the set of cones thus obtained from the congruence classes and let $\F_\Theta$ be the collection of cones consisting of 
arbitrary intersections of the cones in $\C$.
If is convenient to blur the distinction between cones of $\F_\Theta$ and $\Theta$-equivalence classes.
We now proceed to prove Theorem~\ref{main gen} by a series of propositions.
Specifically, Proposition~\ref{fan} verifies that $\F_\Theta$ is a complete fan refined by~$\F$.
Propositions~\ref{fan lat} and~\ref{induced} establish (i) and (ii) respectively.
Assertion (iii) follows from (ii) by Proposition~\ref{induced shelling}.
Proposition~\ref{bisimplicial} proves the claim of the bisimplicial property when~$\A$ is simplicial.
Proposition~\ref{homotopy nonfacial} shows that non-facial intervals are contractible and non-atomic, and 
Proposition~\ref{homotopy facial} completes the proof of (iv) by showing that when~$\A$ is simplicial, facial intervals are atomic and 
homotopy equivalent to spheres of the correct dimensions.

\begin{prop}
\label{fan}
$\F_\Theta$ is a complete fan which is refined by $\F$.
\end{prop}

\begin{proof}
We check the conditions of Lemma~\ref{fan lemma}.
First, suppose that $C_1,C_2\in\C$ intersect in dimension $d-1$.
Then since each is a finite union of regions of~$\A$, there are regions $R_1$ and $R_2$, intersecting in dimension $d-1$ with 
$R_1\subseteq C_1$ and $R_2\subseteq C_2$.
The intersection of $R_1$ and $R_2$ is contained in some hyperplane $H$ of~$\A$, and without loss of generality $R_1\covers R_2$.
Also, $C_1\cap H$ and $C_2\cap H$ are $(d-1)$-dimensional faces of $C_1$ and $C_2$ respectively.
We claim that in fact $C_1\cap H=C_2\cap H$.
To prove this, consider the set $\R$ of~$\A$-regions contained in $C_1$ which intersect $H$ in dimension $d-1$.
Call two regions $Q$ and $R$ in $\R$ {\em adjacent} if $Q\cap R\cap H$ has dimension $d-2$.
Since $C_1\cap H$ is in particular a topological ball, any two regions in $\R$ are connected via a sequence of adjacent elements of $\R$.
Let $Q_1$ be adjacent to $R_1$ in $\R$ and let $Q_2$ be the~$\A$-region whose intersection with $Q_1$ is $Q_1\cap H$, so that in particular
$Q_1\covers Q_2$.
Let $\A'$ be the set of hyperplanes containing $R_1\cap Q_1$.
Then $\A'$ is an arrangement of rank two.
Let $B'$ be the $\A'$-region containing $B$.
If $Q_1$ and $R_1$ are unrelated in $\Po(\A,B)$ then $Q_1\meet R_1$ is some~$\A$-region contained in $B'$, and in particular, 
$Q_1\meet R_1$ is across $H$ from both $Q_1$ and $R_1$.  
But $Q_1\meet R_1$ is congruent to $R_1$, contradicting the fact that $H$ defines a facet of the union over the congruence class of $R_1$.
If $Q_1\le R_1$ then since $S(R_2)=S(R_1)-\set{H}$ and $S(Q_2)=S(Q_1)-\set{H}$, we have $Q_1\meet R_2=Q_2$.
Thus the fact that $Q_1\equiv R_1$ means that $Q_1\meet R_2\equiv R_1\meet R_2$, or in other words, $Q_2\equiv R_2$.
If $Q_1\ge R_1$ we argue similarly that $Q_2\equiv R_2$.
Thus for every region in $\R$, the~$\A$-region whose intersection with $R$ is $R\cap H$ is in $C_2$, so $C_1\cap H\subseteq C_2\cap H$.
By symmetry, we have $C_1\cap H=C_2\cap H$.

Now let $C,D\in\C$   intersect in dimension $k<d-1$.
Then there are~$\A$-regions $Q$ and $R$ with $Q\subseteq C$ and $R\subseteq D$ such that $Q\cap R$ has dimension $k$.
By Lemma~\ref{face}, there is a sequence of $Q=R_0,R_1,\ldots R_t=R$ of regions such that for every $i$ the intersection 
$R_i\cap R_{i-1}$ is $(d-1)$-dimensional and $Q\cap R\subseteq R_i$ for every $i$.
For each $R_i$, let $C_i$ be the cone in $\C$ containing $R_i$.
Then, $C=C_0,C_1,\ldots,C_t=D$ is a sequence of cones in $\C$ such that each cone contains $Q\cap R$ and for each $i$ we have either
$C_i=C_{i-1}$ or $C_i\cap C_{i-1}$ is $(d-1)$-dimensional.

Now we show by induction that $C_0\cap C_1\cap\cdots\cap C_i$ is a face of $C_i$ for each $i$.
The base case is trivial, and when $C_{i-1}\neq C_i$, the intersection $C_{i-1}\cap C_i$ is a $(d-1)$-dimensional face of $C_{i-1}$ and of 
$C_i$, so $C_0\cap C_1\cap\cdots\cap C_i$ is a face of $C_{i-1}\cap C_i$, and in particular a face of $C_i$.
Thus $C_0\cap C_1\cap\cdots\cap C_t$ is a face of $D=C_t$, and since $C\cap D$ has dimension $k$ and each $C_i$ contains the
$k$-dimensional set $Q\cap P$, we have that $C_0\cap C_1\cap\cdots\cap C_k$ is a $k$-dimensional face of $D$.
Because $C$ and $D$ are both convex sets, the intersection $C\cap D$, being $k$-dimensional, cannot be any larger than the $k$ dimensional
face of $D$ it contains, so $C\cap D$ is a face of $D$, and by symmetry, $C\cap D$ is a face of $C$.

We have shown that the intersection of two cones in $\C$ is a face of each.
By Lemma~\ref{fan lemma}, $\F_\Theta$ is a fan.
By construction, $\F_\Theta$ is refined by $\F$.
Since the union of the maximal cones of $\F_\Theta$ is equal to the union of the maximal cones of the complete fan $\F$, the fan $\F_\Theta$ is also complete.
\end{proof}
If $\Theta$ and $\Phi$ are congruences such that $\Theta\le\Phi$ in $\Con(\Po(\A,B))$, then $\F_\Phi$ is refined by $\F_\Theta$.

\begin{prop}
\label{fan lat}
$(\F_\Theta,\Po(\A,B)/\Theta)$ is a fan lattice.
\end{prop}
\begin{proof}
If $[[x]_\Theta,[y]_\Theta]$ is an interval in $\Po(\A,B)/\Theta$, then the union of the corresponding maximal cones of $\F_\Theta$ is equal
to the union of the maximal cones of $\F$ in $[\pidown x,\piup y]$, which is a polytopal cone by the fact that $(\F,\Po(\A,B))$ is a fan lattice.

Each $C\in\C_F$ is the union over an equivalence class of regions in $\Po(\A,B)$, and $F$ is the intersection $\cap_{C\in \C_F}C$.
Since this is an intersection of finite unions, there is some set of representatives $\set{R_C:C\in \C_F}$ whose intersection is full-dimensional in 
$F$.
Since $(\F,\Po(\A,B))$ is a fan lattice, we can let $[R_1,R_2]$ be the interval in $\Po(\A,B)$ consisting of all regions containing 
$\cap_{C\in\C_F}R_C$.
Any $Q\in[R_1,R_2]$ is in some congruence class $C$ containing the full-dimensional subset $\cap_{C\in\C_F}R_C$ of $F$, and thus 
containing $F$.
Thus the set $\C_F$ is the set $\set{[Q]_\Theta:Q\in[R_1,R_2]}$, which by Lemma~\ref{int} is the interval $[[R_1]_\Theta,[R_2]_\Theta]$ 
in $\Po(\A,B)/\Theta$.
\end{proof}

\begin{prop}
\label{induced}
$(\F_\Theta,\Po(\A,B)/\Theta)$ is induced by any linear functional $b$ whose minimum on the unit sphere lies in the interior of $B$.
\end{prop}
\begin{proof}
Suppose $C_1\covered C_2$ in $\Po(\A,B)/\Theta$, and let $R_1\subseteq C_1$ and $R_2\subseteq C_2$ be regions of~$\A$ such that $R_1\cap R_2$ is 
$(d-1)$-dimensional.
Then in particular $R_1\covered R_2$, and since normal vectors to $C_1\cap C_2$ are exactly normal vectors to $R_1\cap R_2$, the result 
follows from the fact that $\Po(\A,B)$ is induced by $b$.
\end{proof}

\begin{prop}
\label{bisimplicial}
If~$\A$ is simplicial, then $(\F_\Theta,\Po(\A,B)/\Theta)$ is bisimplicial with respect to any linear functional $b$ whose minimum on the unit 
sphere lies in the interior of $B$.
\end{prop}
\begin{proof}
Let $C$ be a maximal cone of $\F_\Theta$, and $N_-(C)$ be the set of outward-facing unit normals $\nu$ to facets of $C$ for which 
$b(\nu)<0$.
Then since $(\F_\Theta,\Po(\A,B)/\Theta)$ is induced by $b$, the facets of $C$ corresponding to normals in $N_-(C)$ are the facets 
separating $C$ from maximal cones covered by $C$ in $\Po(\A,B)/\Theta$.
Let $R$ be the region of~$\A$ contained in $C$ such that $\pidown R=R$.
Proposition~\ref{quotient covers} implies that the normals $N_-(C)$ are in one-to-one correspondence with the normals in $N_-(R)$.
Since $R$ is minimal among regions contained in $C$, each facet corresponding to a normal in $N_-(R)$ is contained in a facet 
corresponding to a normal in $N_-(C)$, so that $N_-(C)=N_-(R)$.
Since~$\A$ is simplicial, the set $N_-(C)$ is linearly independent.
The dual argument proves that $N_+(C)$ is linearly independent as well.
\end{proof}
This proof of Proposition~\ref{bisimplicial} goes through under the  weaker hypotheses that $(\F,\Po(\A,B))$ is bisimplicial with respect to any 
linear functional $b$ whose minimum on the unit sphere lies in the interior of $B$.

\begin{prop}
\label{homotopy nonfacial}
If $[x,y]$ is a non-facial interval in $(F_\Theta,\Po(\A,B)/\Theta)$, then $[x,y]$ is a non-atomic interval and $(x,y)$ is contractible.
\end{prop}
\begin{proof}
Let $[x,y]=[[R_1]_\Theta,[R_2]_\Theta]$ be a non-facial interval in $\Po(\A,B)/\Theta$.
Let $R$ be maximal among regions in $[R_1]_\Theta$ which are below $\pidown R_2$.
Thus $[[R_1]_\Theta,[R_2]_\Theta]=[[R]_\Theta,[\pidown R_2]_\Theta]$, so by Lemma~\ref{int}, the interval $[[R_1]_\Theta,[R_2]_\Theta]$
is equal to the set $\set{[Q]_\Theta:Q\in[R,\pidown R_2]}$.
We first show that $[R,\pidown R_2]$ is not a facial interval in $(\F,\Po(\A,B))$.

Suppose for the sake of contradiction that there is a cone $C$ of $\F$ such that $[R,\pidown R_2]$ is exactly the set of maximal cones of $\F$ containing 
$C$.
Let $D$ be minimal among cones of $\F_\Theta$ containing $C$, so that in particular $C$ intersects the relative interior of $D$.
Thus a maximal cone of $\F_\Theta$ contains $D$ if and only if it contains $C$.
If $Q\in[R,\pidown R_2]$ then $[Q]_\Theta$ contains $C$ and therefore also $D$.
Conversely if $[Q]_\Theta$ contains $D$, then $[Q]_\Theta$ contains $C$ and therefore some region $R'\in[Q]_\Theta$ contains $C$,
so that $R'\in[R,\pidown R_2]$.
Thus $\set{[Q]_\Theta:Q\in[R,\pidown R_2]}$ is exactly the set of maximal cones of $\F_\Theta$ containing~$D$, contradicting the fact that 
$[[R_1]_\Theta,[R_2]_\Theta]$ is a non-facial interval.

This contradiction shows that $[R,\pidown R_2]$ is not a facial interval in $\Po(\A,B)$, so that the proper part of $[R,\pidown R_2]$ is 
contractible.
Now since $\pidown R_2$ is minimal in its congruence class and $R$ is maximal in its congruence class among elements $\le\pidown R_2$, the 
restriction of $\Theta$ to $[R,\pidown R_2]$ does not contract any atoms, nor does it set any coatoms equivalent to $\pidown R_2$.
By Lemma~\ref{int}, $[[R_1]_\Theta,[R_2]_\Theta]$ is isomorphic to $[R,\pidown R_2]/\Theta$, and by Corollary~\ref{cross cor} the proper
part of $[R,\pidown R_2]/\Theta$ is contractible as well.
Since $[R,\pidown R_2]$ is not a facial interval in $\Po(\A,B)$, it is not atomic, and since the restriction of $\Theta$ to $[R,\pidown R_2]$ 
neither contracts atoms, nor sets coatoms equivalent to $\pidown R_2$, the interval $[R,\pidown R_2]/\Theta$ is also not atomic.
\end{proof}

\begin{prop}
\label{homotopy facial}
If~$\A$ is simplicial and $[x,y]\subseteq(\F_\Theta,\Po(\A,B)/\Theta)$ is a facial interval associated to a cone $C$ of dimension $k$, then 
$[x,y]$ is an atomic interval with $d-k$ atoms and $(x,y)$ is homotopy equivalent to a $(d-2-k)$-sphere.
\end{prop}
\begin{proof}
As in the proof of Proposition~\ref{fan lat}, we let $[x,y]$ be $[[R_1]_\Theta,[R_2]_\Theta]$, where $[R_1,R_2]$ is a facial interval in 
$\Po(\A,B)$ corresponding to a cone $D$ of $\F$, such that $D$ is a full-dimensional subset of $C$.
Since $[R_1,R_2]$ is a facial interval in $\Po(\A,B)$, $(\F|_C,[R_1,R_2])$ is a simplicial homotopy facial and atomic-facial lattice fan, 
where $\F|_C$ is the fan obtained from $\F$ by restriction.
To obtain the restriction, we have a choice of $p$ in the relative interior of $C$ and a ball $U$ at $p$ intersecting only cones which contain $C$.
We can choose $p$ in the relative interior of $D$, so we have $\F|_C$ refining $(\F_\Theta)|_D$. 
Thus we might as well assume that $[R_1,R_2]=\Po(\A,B)$ and thus that $D=\cap\F$ is a full-dimensional subset of $C=\cap\F_\Theta$.
But since these are both subspaces, we have $C=D$.
We apply Proposition~\ref{may as well} to show that $[B]_\Theta=\set{B}$ and equivalently that $[-B]_\Theta=\set{-B}$.

The atoms of $\Po(\A,B)/\Theta$ number $d-\dim C$ and are in one-to-one correspondence with the atoms of $\Po(\A,B)$.
Since the join of the atoms of $\Po(\A,B)$ is $-B$, the join of the atoms of $\Po(\A,B)/\Theta$ is $\set{-B}$, so $[x,y]$ is atomic and has 
$d-\dim C$ atoms.
Also, by Corollary~\ref{cross cor}, the proper part of $\Po(\A,B)/\Theta$ is homotopy equivalent to the proper part of $\Po(\A,B)$, which is 
homotopy equivalent to a sphere of the desired dimension.
\end{proof}

\begin{prop}
\label{int F}
If $[B]_\Theta=\set{B}$ then $(\cap\F_\Theta)=(\cap\F)$.
\end{prop}
\begin{proof}
For any $\Theta$ the inclusion $(\cap\F)\subseteq(\cap\F_\Theta)$ follows immediately from the fact that each maximal cone of $\F_\Theta$ is
a union of maximal cones of $\F$.
Let $C$ be the maximal cone of $\F_\Theta$ containing $B$.
If $B$ contains $\cap\F_\Theta$ then since $\cap\F_\Theta$ is a subspace, it is contained in every face of $B$, including $\cap\F$, so that 
$(\cap\F)=(\cap\F_\Theta)$.
But $C$ contains $\cap\F_\Theta$.
Thus if $(\cap\F)\subsetneq(\cap\F_\Theta)$ we have $B\subsetneq C$ and thus $[B]_\Theta\neq\set{B}$.
\end{proof}

The proof of Proposition~\ref{bisimplicial} established, independent of the hypothesis that~$\A$ is simplicial, that the facets of the cone 
corresponding to $[R]_\Theta$ in $(\F_\Theta,\Po(\A,B)/\Theta)$ are in one-to-one correspondence with the elements either covered 
by $\pidown R$ or covering $\piup R$ in $\Po(\A,B)$.
Thus we have the following proposition.
Recall that $d$ is the dimension of the vector space in which $\F$ is defined.
\begin{prop}
\label{simp cond}
$\F_\Theta$ is simplicial if and only if for every $[R]_\Theta$ the number of elements covered by $\pidown R$ plus the number of elements 
covering $\piup R$ is $d-\dim(\cap\F_\Theta)$.
\end{prop}
When~$\A$ is simplicial, this condition is equivalent to requiring that the number of elements covered by $\pidown R$ minus the number of 
elements covered by $\piup R$ is $\dim(\cap\F)-\dim(\cap\F_\Theta)$, which is equal to the number of atoms contracted by $\Theta$.

In the case where $\F_\Theta$ is simplicial, Proposition~\ref{shelling} allows the $h$-vector of the corresponding simplicial sphere
to be determined directly from $L/\Theta$.
In particular, the quantity $h_1$ is the number of join-irreducible elements of $\Po(\A,B)/\Theta$.
These are exactly the join-irreducible elements of $\Po(\A,B)$ not contracted by $\Theta$.
When $\Theta$ is a congruence which does not contract any atoms, $f_0$ is the number of atoms plus the number of uncontracted 
join-irreducibles.
By a slight abuse of terminology we call a cone in a complete fan $\F$ a {\em ray} if its dimension is one more than the dimension of the minimal cone
in $\F$.
If the minimal cone of $\F$ is the origin then this is the usual definition of a ray, and if not then we mod out by the minimal cone $\cap\F$ to 
obtain a fan whose minimal cone is the origin.
For a simplicial fan $\F$, the quantity $f_0$ is the number of vertices of the associated simplicial sphere, that is, the number of rays of 
$\F_\Theta$.
Thus the number of rays of $\F_\Theta$ is the number of atoms plus the number of uncontracted join-irreducibles.
By the previous paragraph, the join-irreducibles not contracted by $\Theta$ are in bijection with join-irreducibles $\gamma$ with 
$\piup\gamma=\gamma$.

If we lift the requirement that $\F_\Theta$ be simplicial, when $\Theta$ does not contract any atoms the number of rays of $\F_\Theta$ 
is still the number of atoms of $\Po(\A,B)$ plus the number of join-irreducibles $\gamma$ with $\piup\gamma=\gamma$.
We identify these rays explicitly.
\begin{prop}
\label{rays}
Suppose that~$\A$ is simplicial and that $[B]_\Theta=\set{B}$.
Then the rays of $\F_\Theta$ are exactly the cones arising in one of the following two ways:
\begin{enumerate}
\item[(i) ]For a facet hyperplane $H$ of $B$, let $L$ be the subspace which is the intersection of the other facet hyperplanes of $B$.
Then the cone consisting of points in $L$ weakly separated from $-B$ by $H$ is a ray of $\F_\Theta$.
\item[(ii) ]Given a join-irreducible $\gamma$ of $\Po(\A,B)$ such that $\piup(\gamma)=\gamma$, let $L$ be the intersection of the upper 
facet hyperplanes of $\gamma$.
The cone consisting of points in $L$ weakly separated from $B$ by the unique lower facet hyperplane of $\gamma$ is a ray of $\F(\Theta)$.
\end{enumerate}
\end{prop}
\begin{proof}
Since $[B]_\Theta=\set{B}$, by Proposition~\ref{int F} we have $\cap\F=\cap\F_\Theta$.
By Proposition~\ref{homotopy facial}, if $[x,y]$ is a facial interval in $(\F_\Theta,\Po(\A,B)/\Theta)$ associated to a ray $C$, then $[x,y]$ has 
$d-\dim(\cap\F_\Theta)-1$ atoms and is the quotient modulo $\Theta$ of a facial interval $[x',y']$ with $d-\dim(\cap\F_\Theta)-1$ atoms, 
associated to a ray $D$ in $(\F,\Po(\A,B))$.
Furthermore no atoms of $[x',y']$ are contracted to $x'$ by $\Theta$.
Each $D$ is a half-subspace, that is, the product of $\cap\F$ with a ray (in the usual geometric sense) in $(\cap\F)^\perp$.
Similarly each $C$ is the product of $\cap\F_\Theta$ with a ray in $(\cap\F_\Theta)^\perp$.
Since $D\subseteq C$ and $\cap\F=\cap\F_\Theta$, we have $D=C$.
Thus every interval $[x',y']$ of $\Po(\A,B)$ with $d-\dim(\cap\F)-1$ atoms, none of which are contracted to $x'$, gives rise to a distinct ray of
$\F_\Theta$.

Since~$\A$ is simplicial, the atomic intervals of $\Po(\A,B)$ with $d-\dim(\cap\F)-1$ atoms are of two types.
First, $[B,R]$ where $R$ is the join of a set containing all but one of the atoms of $\Po(\A,B)$, and second, for each join-irreducible 
$\gamma$, the interval $[\gamma,Q]$, where $Q$ is the join of the elements covering $\gamma$.
Since no atoms of $\Po(\A,B)$ are contracted by $\Theta$, no atoms of $[B,R]$ are contracted.
Requiring that no atoms of $[\gamma,Q]$ are contracted to $\gamma$ is exactly the requirement that $\piup\gamma=\gamma$.
It is now easily checked that these rays match the descriptions in (i) and~(ii).
\end{proof}

\section{Weak order on the symmetric group}
\label{symmetric}
For the remainder of this paper we will be concerned with a particular poset of regions, which appears in the guise of the weak order on the 
symmetric group.
In this section we give a brief description of the weak order on the symmetric group, and quote some results concerning its lattice congruences.
Further information, in the more general context of Coxeter groups, can be found in~\cite{Bourbaki,Humphreys} and in 
Section~2 of~\cite{Generalized}.

Let $S_n$ be the symmetric group of permutations of $[n]:=\set{1,2,\ldots,n}$ and write an element $x\in S_n$ in {\em one-line notation}
 $x_1x_2\cdots x_n$, meaning that $x_i:=x(i)$.
The {\em inversion set} $I(x)$ of~$x$ is
\[I(x):=\set{(x_i,x_j):x_i<x_j,i>j}.\]
The length of a permutation~$x$ is $l(x):=|I(x)|$.
Later, we consider permutations in~$S_n$ with $n$ varying.
The {\em size} of a permutation $x$ will denote the $n$ such that $x\in S_n$.

One definition of the right weak order is that $x\le y$ if and only if $I(x)\subseteq I(y)$.
Equivalently, moving up by a cover relation in the right weak order on $S_n$ corresponds to switching adjacent entries in a permutation 
so as to create an inversion.
For the rest of the paper, the phrase ``weak order'' means right weak order, and the symbol ``$S_n$'' denotes the symmetric group as a poset 
under the weak order.
This partial order is the poset of regions of a Coxeter arrangement of type~A, with the inversion set $I$ corresponding to the separating set $S$.
The arrangement is most easily constructed in $\reals^n$, as the set of hyperplanes normal to the vectors $e_i-e_j$ for $1\le j<i\le n$.
The base region $B$ is most conveniently chosen to be the region consisting of points $p=(p_1,p_2,\ldots,p_n)$ with $p_1\le p_2\le\cdots\le p_n$.
The weak order on $S_n$ has a maximal element $w_0:=n(n-1)\cdots 1$.
We denote the identity permutation $12\cdots n$ by $1_n$.

Given a permutation~$x$, say~$x$ has a {\em right descent} at $i$ if $x_i>x_{i+1}$, and say the {\em right descent set} of~$x$ is the 
subset of $[n-1]$ consisting of right descents of~$x$.
The {\em left descent set} of~$x$ is the set consisting of indices $i\in[n-1]$ such that $i+1$ occurs before $i$ in~$x$.
Join-irreducible elements of $S_n$ are permutations with only one right descent.
For any nonempty subset $A\subseteq [n]$, let $A^c:=[n]-A$ and set $m=\min A$ and $M=\max A^c$.
If $\gamma$ is a join-irreducible element of $S_n$ with unique right descent $i$, then $\gamma$ has $\gamma_i>\gamma_{i+1}$ but 
$\gamma_j<\gamma_{j+1}$ for every other $j\in[n-1]$.
Let $A:=\set{\gamma_{i+1},\gamma_{i+2},\ldots,\gamma_n}$.
This is a bijection between join-irreducibles of $S_n$ and nonempty subsets of $[n]$ with $M>m$.
The inverse map takes $A$ to the permutation whose one-line notation consists of the elements of $A^c$ in increasing order followed
by the elements of $A$ in increasing order. 

The weak order on the symmetric group (or more generally on any finite Coxeter group) is a congruence uniform 
lattice~\cite{boundedref,hyperplane}.
In~\cite{congruence}, the poset of irreducibles of $\Con(S_n)$ is determined explicitly as a partial order on the join-irreducibles of $S_n$.

\begin{theorem}\cite[Theorem 8.1]{congruence}
\label{A shard}
The poset $\Irr(\Con(S_n))$ is the transitive closure of the directed graph in which $\gamma_1\to \gamma_2$ if and only if the corresponding subsets $A_1$ 
and $A_2$ satisfy one of the following:
\begin{enumerate}
\item[(i) ]$A_1\cap[1,M_1)=A_2\cap[1,M_1)$ and $M_2>M_1$, or
\item[(ii) ]$A_1\cap(m_1,n]=A_2\cap(m_1,n]$ and $m_2<m_1$.
\end{enumerate}
\end{theorem}

\begin{theorem}\cite[Theorem 8.2]{congruence}
\label{A covers}
Let $\gamma$ and $\gamma'$ be join-irreducibles with corresponding subsets $A$ and $A'$ and let $m$ and $M$ be associated to $A$ as 
described above.
Then $\gamma$ covers $\gamma'$ in $\Irr(\Con(S_n))$ if and only if $A'$ is one of the following:
\begin{eqnarray*}
&A-\set{M+1}&\mbox{for }M<n,\\
&(A-\set{M+1})\cup\set{M}&\mbox{for }M<n,\\
&A\cup\set{m-1}&\mbox{for }1<m, \mbox{ or}\\
&(A\cup\set{m-1})-\set{m}&\mbox{for }1<m.
\end{eqnarray*} 
\end{theorem}

Given a set $K\subseteq[n-1]$, the {\em parabolic subgroup} $(S_n)_K$ of $S_n$ is the subgroup generated by the transpositions 
$\set{(i,i+1):i\in K}$.
Any $x\in S_n$ has a unique factorization $x=x_K\cdot \!\!\phantom{.}^K\! x$ which maximizes $l(x_K)$ subject to the constraints that 
$l(x_K)+l(\!\!\phantom{.}^K\! x)=l(x)$ and that $x_K\in (S_n)_K$.
The set $\!\!\phantom{.}^K\! S_n:=\set{\!\!\phantom{.}^K\! x:x\in S_n}$, called the {\em left quotient} of $S_n$ with respect to $(S_n)_K$,
is a lower interval in weak order, with maximal element $\!\!\phantom{.}^K\! w_0$.
There is an analogous factorization $x=x^K\cdot x_K$, and $(S_n)^K$ is the {\em right quotient}.
A parabolic subgroup $(S_n)_K$ is also a lower interval in the weak order on $S_n$, and the projection $x\mapsto x_K$ is a lattice 
homomorphism.
The corresponding congruence is a parabolic congruence in the sense of Section~\ref{po}.
The parabolic subgroup $(S_{p+q})_{\br{p}}$ for $\br{p}:=[p+q-1]-\set{p}$ is isomorphic to the direct product $S_p\times S_q$, and the map 
from $S_p\times S_q$ is $(u,v)\mapsto u\times v$, where
\[(u\times v)_i=\left\lbrace\begin{array}{ll}u_i&\mbox{if }1\le i\le p\\p+v_{i-p}&\mbox{if }p+1\le i\le p+q\end{array}\right.\]
The upper interval $(S_{p+q})_{\br{p}}\cdot\!\!\phantom{.}^{\br{p}}\! (w_0)$ is also isomorphic to $S_p\times S_q$, and isomorphism is
$(u,v)\mapsto u\ltimes v$, where
\[(u\ltimes v)_i=\left\lbrace\begin{array}{ll}p+v_i&\mbox{if }1\le i\le q\\u_{i-q}&\mbox{if }q+1\le i\le p+q\end{array}\right.\]
We have $u\times v\le u\ltimes v$ in weak order, and the congruence classes of the parabolic congruence associated to ${\br{p}}$ 
are the intervals $[u\times v,u\ltimes v]$.
The join-irreducibles of $(S_{p+q})_{\br{p}}$ are the elements of the form $\gamma\times 1_q$ for $\gamma$ a join-irreducible of $S_p$ and $1_p\times \gamma$
for $\gamma$ a join-irreducible of $S_q$.
The upward projection associated to the parabolic congruence is $w\mapsto w_{\br{p}}\cdot\!\!\phantom{.}^{\br{p}}\! (w_0)$, which 
restricts to an isomorphism from $(S_{p+q})_{\br{p}}$ to $(S_{p+q})_{\br{p}}\cdot\!\!\phantom{.}^{\br{p}}\! (w_0)$.
The following is a specialization of~\cite[Lemma 6.4]{congruence}.
\begin{lemma}
\label{upper}
For a congruence $\Theta$ on $S_{p+q}$, the restriction of $\Theta$ to $(S_{p+q})_{\br{p}}$ corresponds, by the map 
$w_{\br{p}}\mapsto w_{\br{p}}\cdot\!\!\phantom{.}^{\br{p}}\! (w_0)$, to the restriction of $\Theta$ 
to \mbox{$(S_{p+q})_{\br{p}}\cdot\!\!\phantom{.}^{\br{p}}\! (w_0)$}.
\end{lemma}

Define the {\em support} $\supp(x)$ of a permutation~$x$ to be the minimal $K$ such that~$x$ is in $(S_n)_K$, and let the degree of~$x$ 
be $|\supp(x)|$.
The degree of a join-irreducible $\gamma$ in $S_n$ is the magnitude of its unique descent.
That is, if $\gamma_i>\gamma_{i+1}$, then the degree of $\gamma$ is $\gamma_i-\gamma_{i+1}$.
The poset $\Irr(\Con(S_n))$ is dually ranked by the degree.
In a more general context in~\cite{congruence}, it is shown that if $\gamma_1\le \gamma_2$ in $\Irr(\Con(S_n))$, then $\supp(\gamma_2)\subseteq\supp(\gamma_1)$.
A congruence $\Theta$ on $\Po(\A,B)$ is {\em homogeneous of degree $k$} if it is generated by contracting join-irreducibles of degree $k$.

We conclude the section with an observation that allows us to reconstruct a congruence on $S_n$ explicitly from the set of join-irreducibles 
contracted.
Given a permutation $x=x_1x_2\cdots x_n$ with $x_i>x_{i+1}$ define 
\[A(x,i):=\set{x_j:1\le j\le i,x_j>x_i}\cup\set{x_j:i+1\le j\le n,x_j\ge x_{i+1}}.\]
Then $A(x,i)$ has $M=x_i$ and $m=x_{i+1}<M$, so $A(x,i)$ is associated to a join-irreducible which we denote $\lambda(x,i)$.
The permutation $\lambda(x,i)$ consists of all elements of $A^c(x,i)$ in increasing order followed by all elements of $A(x,i)$ in increasing 
order.

\begin{prop}
\label{associated ji}
If $x\covers y$ in the weak order and this covering relation corresponds to transposing $x_i$ and $x_{i+1}$, then a congruence $\Theta$ 
contracts the edge $x\covers y$ if and only if it contracts the join-irreducible $\lambda(x,i)$.
\end{prop}
\begin{proof}
If two entries $a>b$ are inverted in $\lambda(x,i)$, then $b\in A(x,i)$ and $a\in A^c(x,i)$.
If $a<x_{i+1}$ then so is $b$, contradicting the fact that $b\in A(x,i)$.
Therefore, $a=x_j\le x_i$ for some $1\le j\le i$.
Since $a\le x_i$, we have $b\le x_i$ as well, so the fact that $b\in A(x,i)$ implies that $b=x_k$ for some $i+1\le k\le n$.
In particular, $a$ and $b$ are inverted in~$x$ as well and we have shown that $\lambda(x,i)\le x$ in weak order.
The unique element $\lambda_*(x,i)$ covered by $\lambda(x,i)$ is obtained by undoing the inversion $(x_i,x_{i+1})$, so $\lambda_*(x,i)\le y$.
We have $x\meet \lambda(x,i)=\lambda(x,i)$, $y\meet \lambda(x,i)=\lambda_*(x,i)$, $y\join \lambda(x,i)=x$ and $y\join \lambda_*(x,i)=y$.
Applying the definition of lattice congruence to these four equations shows that $\Theta$ contracts the edge $x\covers y$ if and only if it 
contracts the edge $\lambda(x,i)\covers \lambda_*(x,i)$.
\end{proof}

\section{Translational families of congruences}
\label{trans sec}
In this section, we define translational families of congruences and prove Theorem~\ref{subalgebra}, which relates translational families to 
subalgebras of the Malvenuto-Reutenauer algebra.
We also give a combinatorial characterization of translational families in terms of sets of contracted join-irreducibles.

Let $\K$ be a fixed field.
For $n\ge 0$, let $\K[S_n]$ be the vector space over $\K$ spanned by the elements of $S_n$ and let $\K[S_\infty]:=\bigoplus_{n\ge 0}\K[S_n]$.
For $u\in S_p$ and $v\in S_q$, Malvenuto and Reutenauer \cite{MR} defined the {\em shuffle product} $u\productS v\in S_{p+q}$ to be the 
sum of all shuffles of $u$ and $v$.
This is a graded associative product on $\K[S_\infty]$.
Loday and Ronco~\cite{LRorder} pointed out that the shuffle product can be expressed as
\[u\productS v=\sum_{w\in[u\times v,u\ltimes v]}w.\]
and readers not familiar with shuffles may take this as a definition.
The partial order here is the weak order.
In fact, Loday and Ronco used the left weak order, so our product is in fact the dual product used for example in~\cite{Ag-So}.

The product ``$\productS$'' can be rewritten in terms of parabolic subgroups:
\[u\productS v=\sum_{\substack{x\in S_{p+q}\\x_{\br{p}}=u\times v}} x,\]
where $x_{\br{p}}$ refers to the factorization $x=x_{\br{p}}\cdot \!\!\phantom{.}^{\br{p}}\! x$.
Note that $x_{\br{p}}=u\times v$ if and only if $x_{[p-1]}=u$ and $x_{[p+1,p+q-1]}=v$.

For each $n\ge 0$, let $\Theta_n$ be a lattice congruence on $S_n$ and let $(\pidown)_n$ and $(\piup)_n$ 
be the associated downward and upward projections respectively.
As the subscript $n$ is typically given by the context, we refer to all of these projections simply as $\pidown$ and $\piup$.
Let $Z_n^\Theta=S_n/\Theta_n$.
Since $Z_n^\Theta\cong\pidown(S_n)$ we will think of $Z^\Theta_n$ as the subposet $\pidown(S_n)\subseteq S_n$.
Define a graded vector space 
$\K[Z^\Theta_\infty]:=\bigoplus_{n\ge 0}\K[Z^\Theta_n]$.
We often suppress the superscript $\Theta$ and write $\K[Z_\infty]$.
Define a product on $\K[Z_\infty]$ by setting, for $u\in Z_p$ and $v\in Z_q$,
\begin{equation}
\label{prod}
u\productZ v:=\sum_{\substack{x\in Z_{p+q}\\x_{\br{p}}=u\times v}} x,
\end{equation}
that is, we sum over all shuffles of $u$ and $v$ which are the minimal elements of congruence classes of $\Theta_{p+q}$.
Define a map $c:\K[Z_\infty]\to \K[S_\infty]$ by sending each element $x\in Z_n$ to the sum of the elements of the congruence class of~$x$ in 
$S_n$.
The map $c$ is one-to-one and the inverse map $r$, defined on $c(\K[Z_\infty])$, is the map which fixes~$x$ if 
$\pidown x=x$ and maps~$x$ to zero otherwise.
The names $c$ and $r$ for these maps indicate ``class'' and ``representative'' respectively.
We can write the definition of~$\productZ$ concisely as $u\productZ v=r(u\productS v)$ for $u,v\in\K[Z_\infty]$.

For every $p,q\ge 0$, the congruence $\Theta_p\times\Theta_q$ on $S_p\times S_q$ induces a congruence on $(S_{p+q})_{\br{p}}$ via the
map $(u,v)\mapsto u\times v$.
Call the family $\set{\Theta_n}_{n\ge 0}$ of congruences {\em translational} if for every $p,q\ge 0$, this induced congruence on 
$(S_{p+q})_{\br{p}}$ is equal to the restriction of $\Theta_{p+q}$ to $(S_{p+q})_{\br{p}}$.

\bigskip

\noindent {\bf Theorem~ \ref{subalgebra}.}
{\it 
If $\set{\Theta_n}_{n\ge 0}$ is a translational family then the map $c$ embeds $\K[Z^\Theta_\infty]$ as a subalgebra of $\K[S_\infty]$.
}

\bigskip

\begin{proof}
Suppose that $\set{\Theta_n}_{n\ge 0}$ is a translational family of congruences.
The map $c$ respects the vector-space structure and the grading and is one-to-one.
We check that it respects the products of $\K[Z_\infty]$ and $\K[S_\infty]$.

Let $u\in Z_p$ and $v\in Z_q$.
Then 
\[c(u\productZ v)\,\,\,=\,\,\,\sum_{\substack{x\in Z_{p+q}\\x_{\br{p}}=u\times v}}c(x)
\,\,\,=\,\,\,\sum_{\substack{x\in S_{p+q}\\(\pidown x)_{\br{p}}=u\times v}}x.\]
On the other hand,
\[c(u)\productS c(v)\,\,\,=\,\,\,\sum_{\substack{y:\pidown y=u\\z:\pidown z=v}}\,\,\sum_{\substack{x\in S_{p+q}\\x_{\br{p}}=y\times z}}x.\]
Define $I:=\set{x\in S_{p+q}:(\pidown x)_{\br{p}}=u\times v}$.
Define $J$ to be the set of elements $x\in S_{p+q}$ such that, writing $x_{\br{p}}=y\times z$ for some $y\in S_p$ and $z\in S_q$ we have 
$\pidown y= u$ and $\pidown z=v$.
Then $c(u\productZ v)$ is the sum of the elements of $I$ and $c(u)\productS c(v)$ is the sum of the elements of $J$, and we complete the proof by 
showing that $I=J$.

Suppose $x\in I$ and write $x_{\br{p}}=y\times z$ for some $y\in S_p$ and $z\in S_q$.
We have $x\ge x_{\br{p}}$ and therefore $\pidown x\ge \pidown(y\times z)$.
By the order-preserving projection to $(S_{p+q})_{\br{p}}$ this implies $u\times v\ge (\pidown(y\times z))_{\br{p}}$, but since 
$y\times z\in (S_{p+q})_{\br{p}}$, which is a lower interval in $S_{p+q}$, $\pidown(y\times z)\in S_{p+q}$, so 
$(\pidown(y\times z))_{\br{p}}=\pidown(y\times z)$.
Thus $u\times v\ge \pidown(y\times z)$.
On the other hand, $x\ge\pidown x$, so by the order-preserving projection to $(S_{p+q})_{\br{p}}$ we have $y\times z\ge u\times v$ and 
therefore $\pidown(y\times z)\ge\pidown(u\times v)$. 
Since $\set{\Theta_n}$ is a translational family we have 
\[u\times v\ge(\pidown y)\times(\pidown z)\ge(\pidown u)\times(\pidown v)=u\times v,\]
so $\pidown y=u$ and $\pidown z= v$.
Thus $x\in J$ and we have shown that $I\subseteq J$.

Suppose $x\in J$.
Since $x=(y\times z)\cdot\!\!\phantom{.}^{\br{p}}\! x$ we have $(y\times z)\le x\le(y\times z)\cdot\!\!\phantom{.}^{\br{p}}\! w_0$.
We now apply $\pidown$ to the inequality.
Because $\set{\Theta_n}$ is a translational family, $\pidown(y\times z)=(\pidown y)\times(\pidown z)=u\times v$.
Let $(\pidown')$ be the downward projection associated to the restriction of $\Theta_{p+q}$ to 
$(S_{p+q})_{\br{p}}\cdot\!\!\phantom{.}^{\br{p}}\! w_0$.
By Lemma~\ref{upper}, 
\[(\pidown')[(y\times z)\cdot\!\!\phantom{.}^{\br{p}}\! w_0]=\pidown(y\times z)\cdot\!\!\phantom{.}^{\br{p}}\! w_0=(u\times v)\cdot\!\!\phantom{.}^{\br{p}}\! w_0.\]
We have $\pidown[(y\times z)\cdot\!\!\phantom{.}^{\br{p}}\! w_0]\le(\pidown')[(y\times z)\cdot\!\!\phantom{.}^{\br{p}}\! w_0]$, so
\[u\times v\le\pidown x\le\pidown[(y\times z)\cdot\!\!\phantom{.}^{\br{p}}\! w_0]\le(u\times v)\cdot\!\!\phantom{.}^{\br{p}}\! w_0.\]
To this inequality we apply the order-preserving projection down to $(S_{p+q})_{\br{p}}$, thus obtaining 
$u\times v\le(\pidown x)_{\br{p}}\le u\times v$, so $x\in I$.
We have shown that $I=J$.
\end{proof}

There is a more constructive definition of a translational family.
For $k\le n$, $y\in S_k$ and $x\in S_n$, say~$x$ is a {\em translate} of~$y$ if~$x$ is $1_p\times y\times 1_q$ for some $p\ge 0$ and $q\ge 0$.
In this case~$x$ is join-irreducible in $S_n$ if and only if~$y$ is join-irreducible in $S_k$.
Also, since ``$\times$'' is associative and $1_p=1_1\times\cdots\times 1_1$, an arbitrary translation can be obtained as a sequence of translations, each 
of which increases length by~1.
Call~$x$ {\em untranslated} if there is no permutation~$y$ such that~$x$ is a translate of~$y$, or equivalently if $x_1>1$ and $x_n<n$.
For any permutation~$x$ there is a unique untranslated permutation~$y$ such that~$x$ is a translate of~$y$.
Say a permutation $y\in S_k$ has a {\em cliff} at $j$ if $y_j=k$ and $y_{j+1}=1$.
A join-irreducible $\gamma$ in $S_n$ is untranslated if and only if it has a cliff.
This is equivalent to saying the degree of $\gamma$ is $n-1$, which, since $\Irr(\Con(S_n))$ is dually ranked by degree, is equivalent to 
saying that~$\gamma$ is minimal in $\Irr(\Con(S_n))$.
Let $C$ be a set of join-irreducible permutations of various sizes, each of which is untranslated.
(Recall that the size of a permutation $x$ is the $n$ such that $x\in S_n$.) 
For each $n\ge 0$, denote by $\Tr(C)_n$ the smallest congruence on $S_n$ contracting 
every join-irreducible of $S_n$ which is a translate of some element of $C$.
The family of congruences each of which has a single congruence class is $\set{\Tr(21)_n}$, and the family of congruences for which each 
congruence class is a singleton is $\set{\Tr(\emptyset)_n}$.

\begin{prop}
\label{trans}
A family of congruences is translational if and only if it has the form $\set{\Tr(C)_n}_{n\ge 0}$, where $C$ is a set of join-irreducible 
permutations of various sizes, each of which is untranslated.
\end{prop}
\begin{proof}
Let $\gamma$ be a join-irreducible in $S_n$.
Then $\gamma$ is contained in some parabolic subgroup $(S_n)_{\br{p}}$ with $p+q=n$ if and only if either $\gamma=\gamma'\times 1_q$ for some join-irreducible 
$\gamma'$ in $S_p$ or $\gamma=1_p\times \gamma''$ for some join-irreducible $\gamma''\in S_q$.
If $\gamma$ is not in any parabolic subgroup, then in particular it is untranslated.
Since congruences are determined by the set of join-irreducibles they contract, the requirement that $\set{\Theta_n}$ is a translational family 
is equivalent to the requirement that a join-irreducible is contracted if and only if all of its translates are.
Therefore a translational family is ${\Tr(C)_n}$, where $C$ is the set of untranslated join-irreducibles of various sizes contracted by the family.
\end{proof}
The proof of Proposition~\ref{trans} constructs $C$ as the set of all contracted untranslated join-irreducibles.
However, in many examples we take $C$ to be a finite generating set.
The following lemmas are easily checked by reducing to the case $p+q=1$ and applying Theorem~\ref{A shard}.
\begin{lemma}
\label{lift}
Let $\gamma_1$ and $\gamma_2$ be join-irreducibles in $S_k$.
Then $\gamma_1\le\gamma_2$ in $\Irr(\Con(S_k))$ if and only if $1_p\times\gamma_1\times 1_q\le 1_p\times\gamma_2\times 1_q$ in
$\Irr(\Con(S_{p+k+q}))$.
\end{lemma}
\begin{lemma}
\label{lift 2}
If $\gamma_1\ge 1_p\times\gamma_2\times 1_q$ in $\Irr(\Con(S_n))$ then $\gamma_1=1_p\times\gamma_1'\times 1_q$ for some $\gamma_1'$.
\end{lemma}

Refinements of congruences give rise to further subalgebra relationships.
Specifically, let $\set{\Theta_n}$ and $\set{\Phi_n}$ be two translational families such that $\Phi_n$ refines $\Theta_n$ for each $n$.
Alternately, we can think of $\Theta_n$ as a congruence on the lattice $S_n/\Phi_n$.
Then $\K[Z^\Theta_\infty]$ is a subalgebra of $\K[Z^\Phi_\infty]$.

This restriction of the refinement order on families to translational families is a distributive lattice.
Specifically, the join of two translational families $\set{\Tr(C_1)}$ and $\set{\Tr(C_2)}$ is $\set{\Tr(C_1\cup C_2)}$.
If $C_1$ is the complete set of untranslated join-irreducibles contracted by $\set{\Tr(C_1)}$ and similarly for $C_2$ then
$\set{\Tr(C_1\cap C_2)}$ is the meet of $\set{\Tr(C_1)}$ and $\set{\Tr(C_2)}$.
We wish to define a partial order $\Tr_\infty$ on untranslated join-irreducibles such that the possible sets $C$ of all contracted untranslated 
join-irreducibles for a translational family are exactly the order ideals in $\Tr_\infty$.
A priori, this means defining $\Tr_\infty$ as the transitive closure of the relation setting $\gamma_1\ge\gamma_2$ whenever some translate of 
$\gamma_1$ is above some translate of $\gamma_2$ in $\Irr(\Con(S_n))$ for some $n$.
However, this definition can be simplified.
\begin{prop}
\label{simp def}
Let $\gamma_1$ and $\gamma_2$ be untranslated join-irreducibles, such that $\gamma_2\in S_k$.
Then $\gamma_1\covers\gamma_2$ in $\Tr_\infty$ if and only if some translate $\gamma_1'$ of $\gamma_1$ covers $\gamma_2$ in 
$\Irr(\Con(S_k))$.
Furthermore, $\gamma_1'$ is either $1_1\times\gamma_1$ or $\gamma_1\times 1_1$.
\end{prop}
\begin{proof}
Suppose $\gamma_1\covers\gamma_2$ in $\Tr_\infty$, so that in particular some translate $\gamma_1'$ of $\gamma_1$ is above 
$1_p\times\gamma_2\times 1_q$ in $\Irr(\Con(S_{p+k+q}))$.
Then by Lemma~\ref{lift 2}, $\gamma_1'=1_p\times\gamma_1''\times 1_q$ for some $\gamma_1''$.
By Lemma~\ref{lift} we have $\gamma_1''>\gamma_2$ in $\Irr(\Con(S_k))$.
But $\gamma_1$ is untranslated, and $\gamma''_1$ is not minimal in $\Irr(\Con(S_k))$, so it is a translate.
Thus $\gamma_1\neq\gamma_1''$, so $\gamma_1''$ is a translate of $\gamma_1$.
If there is some $\gamma_3$ such that $\gamma_1''>\gamma_3>\gamma_2$ in $\Irr(\Con(S_k))$, then there is some untranslated 
join-irreducible $\gamma_3'$ such that $\gamma_3$ is a translate of $\gamma_3'$, and $\gamma_1>\gamma_3'>\gamma_2$ in 
$\Tr_\infty$, contradicting the hypothesis that $\gamma_1\covers\gamma_2$.
Thus $\gamma_1''\covers\gamma_2$ in $\Irr(\Con(S_k))$.

Suppose conversely that some translate $\gamma'_1$ covers $\gamma_2$ in $\Irr(\Con(S_n))$.
Recall that $\Irr(\Con(S_n))$ is dually ranked by degree, and that the degree of a join-irreducible is the magnitude of its unique descent.
Since $\gamma_2$ has degree $k-1$, the translate $\gamma_1'$ has degree $k-2$.
Thus the unique descent of $\gamma_1'$ consists either of the entry $k$ followed by $2$ or $k-1$ followed by $1$, so that $\gamma_1'$ is 
either $1_1\times\gamma_1$ or $\gamma_1\times 1_1$.
We have $\gamma_1>\gamma_2$ in $\Tr_\infty$.
If there is some $\gamma_3$ such that $\gamma_1\covers\gamma_3>\gamma_2$ in $\Tr_\infty$, then by the previous paragraph there would 
have to be a translate of $\gamma_1$ covering $\gamma_3$ and a translate of $\gamma_3$ greater than $\gamma_2$.
But this is impossible since $\gamma_1\in S_{k-1}$ and $\gamma_2\in S_k$.
\end{proof}

To explicitly describe the cover relations in $\Tr_\infty$ we introduce an operation called {\em insertion}.
Let $\gamma$ be a join-irreducible in $S_n$ with associated subset $A$ and let $i\in[n+1]$.
Then the {\em left insertion} of $i$ in $\gamma$ is a join-irreducible $\LI_i(\gamma)$ in $S_{n+1}$ whose associated subset is 
$(A\cap[1,i-1])\cup\set{j+1:j\in A\cap[i,n]}$.
The {\em right insertion} of $i$ in $\gamma$ is a join-irreducible $\RI_i(\gamma)$ in $S_{n+1}$ whose associated subset is 
$(A\cap[1,i-1])\cup\set{i}\cup\set{j+1:j\in A\cap[i,n]}$.
When $\gamma$ is written $\gamma_1\gamma_2\cdots\gamma_n$, it consists of the elements of $A^c$ in increasing order on the left, followed
by the elements of $A$ in increasing order on the right.
The effect of these insertions is to increase each entry $\ge i$ by 1 and then insert $i$ into either the left increasing sequence or the right 
increasing sequence.
Note that $\LI_1(\gamma)=1_1\times\gamma$ and $\RI_{n+1}(\gamma)=\gamma\times1_1$.

\begin{figure}[ht]
\caption{The top four ranks of $\Tr_\infty$.}
\label{Tr inf}
\centerline{\epsfbox{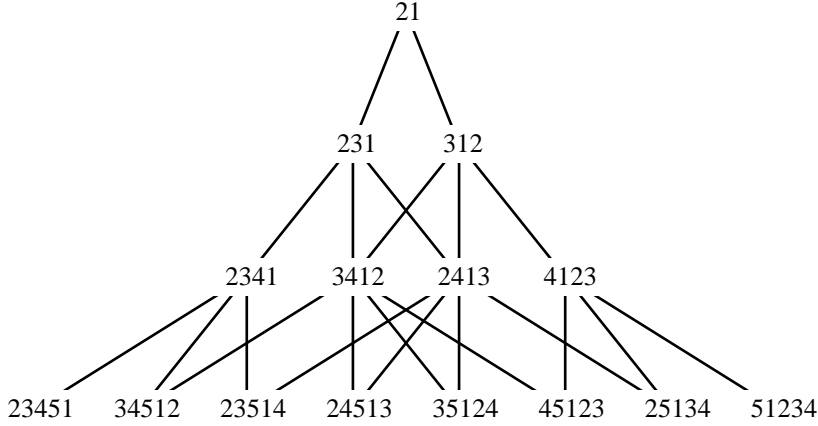}}   
\end{figure}

\begin{prop}
\label{tr}
Let $\gamma$ be an untranslated join-irreducible in $S_n$.
Then the set of elements covered by $\gamma$ in $\Tr_\infty$ is $\set{R_1(\gamma),\LI_2(\gamma),\RI_n(\gamma),\LI_{n+1}(\gamma)}$.
These elements are not necessarily distinct.
\end{prop}
\begin{proof}
We apply Proposition~\ref{simp def}.
Using Theorem~\ref{A covers}, it is easily checked that the two elements covered by $1_1\times\gamma$ are $R_1(\gamma)$ and 
$L_2(\gamma)$, and that the two elements covered by $\gamma\times 1_1$ are $R_n(\gamma)$ and $L_{n+1}(\gamma)$.
\end{proof}
The poset $\Tr_\infty$ is dually ranked by size.
The top four ranks of $\Tr_\infty$ are pictured in Figure~\ref{Tr inf}.
Reflecting this picture through a vertical line is the symmetry of $\Tr_\infty$ which corresponds to applying to each $\Irr(\Con(S_n))$ the 
antipodal symmetry defined in~\cite[Section 6]{congruence}.

\section{Insertional families of congruences}
\label{ins}
In this section we define insertional families of congruences and prove Theorem~\ref{subcoalgebra}, which relates insertional families to 
subcoalgebras of the Malvenuto-Reutenauer Hopf algebra.
We also give a combinatorial characterization of insertional families in terms of sets of contracted join-irreducibles.

Let $\K[X_\infty]$ be a vector space graded by the non-negative integers, such that $X_n$ is a basis for the $n$-graded component.
A {\em (graded) coproduct} on $\K[X_\infty]$ is a map $\Delta:\K[X_\infty]\to \K[X_\infty]\otimes \K[X_\infty]$ 
such that the image of the restriction of $\Delta$ to $\K[X_n]$ is contained in $\bigoplus_{p+q=n}\K[X_p]\otimes \K[X_q]$.
The coproduct $\Delta$ is {\em coassociative} if $(\Delta\otimes I)\circ\Delta=(I\otimes\Delta)\circ\Delta$, where $I$ is the identity.
If $|X_0|=1$ and the orthogonal projection $\ep:\K[X_\infty]\to\K[X_0]\cong\K$ satisfies $(\ep\otimes I)\circ\Delta=1\otimes I$ and 
$(I\otimes \ep)\circ\Delta=I\otimes 1$, where ``$1$'' is the map whose image is $\set{1}\subset\K$, then $\K[X_\infty]$ is called a {\em 
graded, connected coalgebra}.
A {\em subcoalgebra} $\K[Y_\infty]$ of $\K[X_\infty]$ is a graded subspace such that the image of the restriction $\Delta_Y$ of $\Delta$ 
to $\K[Y_\infty]$ is contained in $\K[Y_\infty]\otimes \K[Y_\infty]$.

Let $\K[X_\infty]$ be a graded, connected coalgebra with coproduct $\Delta$.
There is a standard construction of a coproduct on $\K[X_\infty]\otimes\K[X_\infty]$, induced by $\Delta$, which makes 
$\K[X_\infty]\otimes\K[X_\infty]$ a graded, connected coalgebra. 
If $\K[X_\infty]$ is also an algebra whose product ``$\product$'' is a coalgebra homomorphism $\K[X_\infty]\otimes\K[X_\infty]\to\K[X_\infty]$ then $(\K[X_\infty],\product,\Delta)$ is called a {\em graded, connected Hopf algebra}.
For the rest of the paper, the term ``Hopf algebra'' will mean a graded, connected Hopf algebra.

Suppose $(\K[X_\infty],\product,\Delta)$ is a Hopf algebra and $\K[Y_\infty]$ is both a subcoalgebra and a subalgebra of 
$\K[X_\infty]$, with $|X_0|=|Y_0|=1$.
Then $(\K[Y_\infty],\productY,\Delta_Y)$ is a Hopf algebra, and in particular a {\em graded sub Hopf algebra} of 
$(\K[X_\infty],\product,\Delta)$, where ``$\productY$'' is the restriction of ``$\product$'' to $\K[Y_\infty]$.
For more information on Hopf algebras, see~\cite{MM,Sweedler}. 

Malvenuto and Reutenauer \cite{MR} defined a coproduct $\Delta_S$ on $\K[S_\infty]$ such that $(\K[S_\infty],\!\productS\!,\Delta_S)$ 
is a Hopf algebra.
To any sequence $(a_1,a_2,\ldots,a_p)$ of distinct integers, we associate a {\em standard permutation} $\st(a_1,a_2,\ldots,a_p)$.
This is the permutation $u\in S_p$ such that for each $i,j\in[p]$ we have $u_i<u_j$ if and only if $a_i<a_j$.
So for example $\st(73591)=42351$.
The standard permutation of the empty sequence is the empty permutation in $S_0$.
The coproduct of an element $x\in S_n$ is
\[\Delta_S(x)=\sum_{p=0}^n \st(x_1,\ldots,x_p)\otimes\st(x_{p+1},\ldots,x_n).\]
For convenience, let $(S_n)_{\br{0}}$ and $(S_n)_{\br{n}}$ both denote $S_n$ which is, in a trivial sense, a parabolic subgroup of itself.
We should think of $(S_n)_{\br{0}}$ as $S_0\times S_n$ and $(S_n)_{\br{n}}$ as $S_n\times S_0$.
As discussed in Section~\ref{symmetric}, for $p\in[n-1]$, the product $S_p\times S_{n-p}$ is isomorphic to the parabolic subgroup 
$(S_n)_{\br{p}}$ by the map $(u,v)\mapsto u\times v$.
The following formula for $\Delta_S$ is useful in the proof of Theorem~\ref{subcoalgebra} despite the fact that the second and third sums are 
each sums of a single term.
\[\Delta_S(x)=\sum_{p=0}^n\sum_{w\in S_n^{\br{p}}}\,\,\sum_{\substack{u\in S_p,v\in S_{n-p}\\w\cdot(u\times v)=x}}u\otimes v.\]
To avoid confusion, we stress the fact that $S_n^{\br{p}}$ is the right quotient of $S_n$ with respect to $(S_n)_{\br{p}}$ rather than the
left quotient which was used to define the product~$\productS$.

For $\set{\Theta_n}_{n\ge 0}$ with $\pidown$ and $\piup$ be as before, we define a coproduct $\Delta_Z$ on $\K[Z^\Theta_\infty]$.
Given a sequence $x_1,x_2,\ldots,x_p$ of distinct positive integers, let 
$\st^\Theta(x_1,x_2,\ldots,x_p)=\st(x_1,x_2,\ldots,x_p)\in\K[Z^\Theta_\infty]$ if $\pidown(\st(x_1,x_2,\ldots,x_p))=\st(x_1,x_2,\ldots,x_p)$, and
otherwise let $\st^\Theta(x_1,x_2,\ldots,x_p)=0$.
For $x\in Z_n$, let $\Delta_Z=(r\otimes r)\circ\Delta_S\circ c$, where $c$ and~$r$ are the maps defined in Section~\ref{trans sec}.
That is:
\begin{equation}
\label{coprod}
\Delta_Z(x):=\sum_{\substack{y\in S_n\\\pidown y=x}}\,\,\sum_{p=0}^n\st^\Theta(x_1,\ldots,x_p)\otimes\st^\Theta(x_{p+1},\ldots,x_n).
\end{equation}
We can rearrange the sum to read
\[\Delta_Z(x)=\sum_{p=0}^n\,\,\sum_{w\in S_n^{\br{p}}}\,\,\sum_{\substack{u\in Z_p,v\in Z_{n-p}\\\pidown(w\cdot(u\times v))=x}} 
u\otimes v.\]

Each left coset of a parabolic subgroup in $S_n$ is an interval in the weak order isomorphic to the weak order on the parabolic subgroup.
For any $p,q\ge 0$ and any $w\in S_{p+q}^{\br{p}}$, the product $S_p\times S_q$ is isomorphic to $w\cdot(S_{p+q})_{\br{p}}$ by the 
map $(u,v)\mapsto w\cdot(u\times v)$.
The congruence $\Theta_p\times\Theta_q$ on $S_p\times S_q$ induces a congruence on $w\cdot(S_{p+q})_{\br{p}}$ via the map 
$(u,v)\mapsto w\cdot(u\times v)$.
Call the family $\set{\Theta_n}_{n\ge 0}$ of congruences {\em insertional} if for every $p,q\ge 0$ and every $w\in S_{p+q}^{\br{p}}$, the 
congruence induced on $w\cdot(S_{p+q})_{\br{p}}$ by $\Theta_p\times\Theta_q$ is a refinement of the restriction of $\Theta_{p+q}$ to 
$w\cdot(S_{p+q})_{\br{p}}$.
The term ``insertional'' will be justified later by Proposition~\ref{insert def}.

\bigskip

\noindent {\bf Theorem~ \ref{subcoalgebra}.}
{\it 
If $\set{\Theta_n}_{n\ge 0}$ is an insertional family then the map $c$ embeds $\K[Z^\Theta_\infty]$ as a subcoalgebra of $\K[S_\infty]$.
}

\bigskip

\begin{proof}
We check that $c$ respects the coproducts.
Let $x\in Z_n$ and think of~$x$ as usual as an element of $S_n$ with $\pidown x=x$.
Then
\[(c\otimes c)(\Delta_Z(x))=\sum_{p=0}^n\,\,\sum_{w\in S_n^{\br{p}}}\,\,\sum_{\substack{u\in Z_p,v\in Z_{n-p}\\w\cdot(u\times v)\in[x,\piup x]}} c(u)\otimes c(v).\]
On the other hand,
\[\Delta_S(c(x))=\sum_{p=0}^n\,\,\sum_{w\in S_n^{\br{p}}}\,\,\sum_{\substack{u\in S_p,v\in S_{n-p}\\w\cdot(u\times v)\in[x,\piup x]}} u\otimes v.\]
For each $w\in S_n^{\br{p}}$, the map $(u,v)\mapsto w\cdot(u\times v)$ maps each $(\Theta_p\times\Theta_q)$-class in $S_p\times S_q$
to an interval in \mbox{$w\cdot(S_{p+q})_{\br{p}}$}.
If  $\set{\Theta_n}$ is an insertional family, then for each $w$ this interval is either entirely contained in 
$w\cdot(S_{p+q})_{\br{p}}\cap[x,\piup x]$ or disjoint from $w\cdot(S_{p+q})_{\br{p}}\cap[x,\piup x]$.
Thus these two sums are equal.
\end{proof}

\begin{prop}
\label{insert def}
A family of congruences is insertional if and only if for every join-irreducible $\gamma$ contracted by $\Theta_n$ with associated subset $A$ 
and $m=\min A$, $M=\max A^c$, the following two conditions hold:
\begin{enumerate}
\item[(i) ]$\Theta_{n+1}$ contracts the right insertion $\RI_i(\gamma)$ for every $i\in[m+1,M+1]$.
\item[(ii) ]$\Theta_{n+1}$ contracts the left insertion $\LI_i(\gamma)$ for every $i\in[m,M]$.
\end{enumerate}
\end{prop}
\begin{proof}
Let $\gamma$ be a join-irreducible in $S_n$.
Let $A$ be the subset corresponding to $\gamma$ and let the unique descent of $\gamma$ be between the entries $\gamma_k=M$ and 
$\gamma_{k+1}=m$.
The choice of $w\in S_{n+1}^{\br{n}}$ amounts to choosing an element $i\in[n+1]$ so that the one-line notation for $w$ consists of the entries
of $[n+1]-\set{i}$ in increasing order followed by the entry~$i$.
Choose $i\in[m+1,M+1]$ and let~$x$ be the permutation $w\cdot(\gamma\times 1_1)$, so that $\st(x_1,x_2,\ldots,x_n)=\gamma$.
Let $y=w\cdot(\gamma_*\times 1_1)$, so that~$y$ agrees with~$x$ except that the entries in positions $k$ and $k+1$ are transposed.
We have
\[x_j=\left\lbrace\begin{array}{ll}
\gamma_j&\mbox{ if }j<n+1\mbox{ and }\gamma_j<i,\\
\gamma_j+1&\mbox{ if }j<n+1\mbox{ and }\gamma_j\ge i,\mbox{ or}\\
i&\mbox{ if }j=n+1.
\end{array}\right.\]
Consider the join-irreducible $\lambda(x,k)$ defined in connection with Proposition~\ref{associated ji}.
Since $i>m$, $\lambda(x,k)$ is constructed from~$x$ by moving the entry $i$ into a position to the right of position $k$
such that the entries in positions $k+1$ to $n+1$ are increasing.
In other words, $\lambda(x,k)$ is $\RI_i(\gamma)$.
By Proposition~\ref{associated ji}, the edge \mbox{$w\cdot(\gamma\times 1_1)\covers w\cdot(\gamma_*\times 1_1)$} is contracted by $\Theta_{n+1}$
if and only if the join-irreducible $R_i(\gamma)$ is contracted by~$\Theta_{n+1}$.

Choosing $w\in S_{n+1}^{\br{1}}$ amounts to choosing an $i\in[n+1]$ so that the one-line notation for $w$ consists of the entry $i$ followed
by the entries in $[n+1]-\set{i}$ in increasing order.
Choose some $i\in[m,M]$, let $x=w\cdot(1_1\times\gamma)$ and let $y=w\cdot(1_1\times\gamma_*)$.
We have
\[x_j=\left\lbrace\begin{array}{ll}
i&\mbox{ if }j=1,\\
\gamma_{j-1}&\mbox{ if }j>1\mbox{ and }\gamma_j<i,\mbox{ or}\\
\gamma_{j-1}+1&\mbox{ if }j>1\mbox{ and }\gamma_j\ge i.
\end{array}\right.\]
Since $i\le M$, we have $\LI_i(\gamma)=\lambda(x,k)$ which, by Proposition~\ref{associated ji} is contracted if and only
if the edge $w\cdot(1_1\times\gamma)\covers w\cdot(1_1\times\gamma_*)$ is contracted.

If $\set{\Theta_n}$ is insertional then for any join-irreducible $\gamma\in S_n$ contracted by $\Theta_n$, the edges 
$w\cdot(\gamma\times 1_1)\covers w\cdot(\gamma_*\times 1_1)$ and $w\cdot(1_1\times\gamma)\covers w\cdot(1_1\times\gamma_*)$
are contracted by $\Theta_{n+1}$.
By the previous two paragraphs this implies (i) and (ii).

Since congruences are determined by the set of join-irreducibles they contract, the definition of an insertional family can be rewritten as the 
following requirements for each $p,q\ge 0$ and $w\in S_{p+q}^{\br{p}}$.
\begin{enumerate}
\item[(i') ] If a join-irreducible $\gamma\in S_p$ is contracted by $\Theta_p$, then $\Theta_{p+q}$ contracts the edge 
\mbox{$w\cdot(\gamma\times 1_q)\covers w\cdot(\gamma_*\times 1_q)$}, and
\item[(ii') ] If a join-irreducible $\gamma\in S_q$ is contracted by $\Theta_q$ then $\Theta_{p+q}$ contracts the edge 
\mbox{$w\cdot(1_p\times\gamma)\covers w\cdot(1_p\times\gamma_*)$}.
\end{enumerate}
Now suppose conditions (i) and (ii) in the statement of the proposition hold and let $\gamma$ be a join-irreducible in $S_p$ contracted by 
$\Theta_p$ with $A$, $m$, $M$ and $k$ as above.
For fixed $p$ and $q$, the choice of $w\in S_{p+q}^{\br{p}}$ corresponds to choosing some subset $Q$ of $[p+q]$ with $q$ elements.
Let $a_1a_2\cdots a_p$ be the unique permutation of the set $[p+q]-Q$ whose standard permutation is $\gamma$.
Let $x:=w\cdot(\gamma\times 1_q)$, so that~$x$ consists of the entries $a_1\cdots a_p$ followed by the elements of $Q$ in increasing order.
Let $b_1,\ldots,b_r$ be the elements of $Q\cap[1,a_{k+1}]$ in increasing order, let $c_1,\ldots,c_s$ be the elements of 
$Q\cap[a_{k+1}+1,a_k+1]$ in increasing order, and let $d_1,\ldots,d_t$ be the elements of $Q\cap[a_k+2,n+1]$ in increasing order.
Then 
\[\lambda(x,k)=\RI_{d_t}\cdots\RI_{d_1}(\RI_{c_s}\cdots\RI_{c_1}(\LI_{b_k}\cdots\LI_{b_1}(\gamma))).\]
For any join-irreducible $\gamma$ with corresponding $m$ and $M$, if $j<m$ we have $\LI_j(\gamma)=\LI_m(\gamma)$ so that, by condition 
(ii), $\LI_{b_k}\cdots\LI_{b_1}(\gamma)$ is contracted.
By condition (i), $\RI_{c_s}\cdots\RI_{c_1}(\LI_{b_k}\cdots\LI_{b_1}(\gamma))$ is contracted.
For any join-irreducible $\gamma$ with corresponding $m$ and $M$, if $j>M+1$ we have $\RI_j(\gamma)=\RI_{M+1}(\gamma)$, so by 
condition (ii), $\lambda(x,k)$ is contracted.
Thus by Proposition~\ref{associated ji}, $w\cdot(\gamma\times 1_q)\covers w\cdot(\gamma_*\times 1_q)$ is contracted.
We have verified that conditions (i) and (ii) imply condition (i').
The proof for (ii') is similar.
\end{proof}

\section{Sub Hopf algebras}
\label{H}
Recall that an {\em $\scrH$-family} is a translational and insertional family of congruences, and that Corollary~\ref{subHopf} states that 
if $\set{\Theta_n}_{n\ge 0}$ is an $\scrH$-family then the map $c$ embeds $\K[Z^\Theta_\infty]$ as a sub Hopf algebra of $\K[S_\infty]$.
In this section we study $\scrH$-families.

\begin{prop}
\label{trins}
Let $C$ be a set of untranslated join-irreducibles of various sizes.
Then $C$ is the complete set of untranslated join-irreducibles contracted by an $\scrH$-family if and only $C$ is closed under insertions 
which are not translations.
\end{prop}
\noindent In other words, for every $\gamma\in C\cap S_n$, the requirement is that $\set{\RI_i(\gamma):i\in[1,n]}\subset C$ and 
$\set{\LI_i(\gamma):i\in[2,n+1]}\subset C$.
\begin{proof}
By Proposition~\ref{tr}, $C$ is the complete set of untranslated join-irreducibles contracted by a translational family if and only if for every 
$\gamma\in C\cap S_n$ we have $\set{R_1(\gamma),\LI_2(\gamma),\RI_n(\gamma),\LI_{n+1}(\gamma)}\subset C$.
If the family is insertional as well, by Proposition~\ref{insert def}, the join-irreducibles $\RI_i(\gamma)$ for $i\in[2,n+1]$ and 
$\LI_i(\gamma)$ for $i\in[1,n]$ are contracted as well.
Note that $\RI_{n+1}(\gamma)=\gamma\times 1_1$ and $\LI_1(\gamma)=1_1\times\gamma$, but that all of these other insertions yield
untranslated join-irreducibles.
Thus $\set{\RI_i(\gamma):i\in[1,n]}\subseteq C$ and $\set{\LI_i(\gamma):i\in[2,n+1]}\subseteq C$.

Conversely, let $C$ have $\set{\RI_i(\gamma):i\in[1,n]}\subseteq C$ and $\set{\LI_i(\gamma):i\in[2,n+1]}\subseteq C$ for every $n$ 
and every $\gamma\in C\cap S_n$.
Let $\gamma$ be a join-irreducible contracted by $\Theta_n$ with associated subset $A$, let $m=\min A$ and let $M=\max A^c$.
Then $\gamma=1_{m-1}\times\gamma'\times 1_{n-M}$ for some untranslated join-irreducible $\gamma'$.
If $i\in[m,M+1]$, then $\RI_i(\gamma)=1_{m-1}\times\RI_{i-m+1}(\gamma')\times 1_{n-M}$, and 
$\LI_i(\gamma)=1_{m-1}\times\LI_{i-m+1}(\gamma')\times 1_{n-M}$, so these conditions on untranslated join-irreducibles imply 
conditions (i) and (ii) of Proposition~\ref{insert def}.
\end{proof}

We define a partial order $\scrH_\infty$, in analogy to $\Tr_\infty$, such that the possible sets $C$ of contracted untranslated 
join-irreducibles for an $\scrH$-family are exactly the order ideals in $\scrH_\infty$.
In particular, the set of $\scrH$-families is a distributive lattice.
Proposition~\ref{trins} can be interpreted as a description of the cover relations in $\scrH_\infty$, keeping in mind that the set
$\set{\RI_i(\gamma):i\in[1,n]}\cup\set{\LI_i(\gamma):i\in[2,n+1]}$ does not necessarily contain $2n$ distinct elements.
The poset $\scrH_\infty$ is an extension of $\Tr_\infty$, in the sense that the underlying sets coincide and every order relation in $\Tr_\infty$
is an order relation in $\scrH_\infty$.
Like $\Tr_\infty$, the poset $\scrH_\infty$ is dually ranked by size.
In Figure~\ref{Tr inf}, one adds in the cover relations $2341\covers 24513$ and $4123\covers 35124$ to obtain a diagram of the top four 
ranks of $\scrH_\infty$.

Given permutations $y=y_1y_2\cdots y_k\in S_k$ and $x=x_1x_2\cdots x_n\in S_n$, say that the {\em pattern}~$y$ 
{\em occurs} in~$x$ if there are integers $1\le i_1<i_2<\cdots<i_k\le n$ such that  for all $1\le p<q\le k$ we have $y_p<y_q$ 
if and only if $x_{i_p}<x_{i_q}$.
Otherwise, say that~$x$ {\em avoids}~$y$.
For more information on patterns in permutations, see \cite{Wilf}.
We extend the definition of pattern avoidance in order to describe $\scrH$-families.
Recall that a permutation $y\in S_k$ has a cliff at $j$ if $y_j=k$ and $y_{j+1}=1$.
If~$y$ has a cliff at $j$, say~$y$ occurs in~$x$ {\em with adjacent cliff} if there is some occurrence $x_{i_1}x_{i_2}\cdots x_{i_k}$ of~$y$ in~$x$
such that $i_{j+1}=i_j+1$. 
Otherwise say~$x$ {\em avoids cliff-adjacent instances} of~$y$.
For an untranslated join-irreducible $\gamma\in S_k$ with a cliff at $j$, that is $\gamma_j=k$ and $\gamma_{j+1}=1$, a {\em scramble} of 
$\gamma$ is any permutation $\sigma$ with $\sigma_j=k$, $\sigma_{j+1}=1$ and $\set{\sigma_i:i\in[j-1]}=\set{\gamma_i:i\in[j-1]}$.
Notice that every scramble of $\gamma$ has a cliff at~$j$.

Let $C$ be a set of join-irreducible permutations of various sizes, each of which is untranslated, and define $\set{\scrH(C)_n}$ to be the
smallest $\scrH$-family of congruences contracting each element of $C$.
Thus the complete set of untranslated join-irreducibles contracted by $\set{\scrH(C)_n}$ is the smallest order ideal of $\scrH_\infty$ 
containing $C$.

\begin{prop}
\label{avoid}
A join-irreducible $\gamma'\in S_n$ is contracted by $\scrH(C)_n$ if and only if there is some $\gamma\in C$ which occurs as a pattern
in $\gamma'$.
\end{prop}
\begin{proof}
Since $\scrH$-families correspond to order ideals in $\scrH_\infty$, we may as well take $C=\set{\gamma}$ for some untranslated 
join-irreducible $\gamma\in S_k$.
Other order ideals are obtained as unions of these principal order ideals.
Also, we can reduce to the case where $\gamma'$ is untranslated.
Otherwise, write $\gamma'=1_p\times\gamma''\times 1_q$ for some $p$ and $q$ and some untranslated join-irreducible $\gamma''$.
Because $\set{\scrH(C)_n}$ is in particular a translational family, $\gamma'$ is contracted if and only if $\gamma''$ is contracted.
Furthermore, $\gamma'$ contains the untranslated join-irreducible $\gamma$ if and only if $\gamma''$ contains $\gamma$.
Now it is easily proven by induction on $n$ that $\gamma'$ is contracted by $\scrH(\gamma)_n$ if and only if $\gamma$ occurs in $\gamma'$. 
\end{proof}

\begin{theorem}
\label{c-a avoid}
A permutation $x\in S_n$ is contracted by $\scrH(C)_n$ if and only for some $\gamma\in C$ there is a scramble $\sigma$ of $\gamma$ 
which occurs in~$x$ with adjacent cliff.
\end{theorem}
\begin{proof}
Again we reduce to the case where $C=\set{\gamma}$ for some $\gamma\in S_k$, because when $|C|>1$, the permutation~$x$ is 
contracted by $\scrH(C)_n$ if and only if it is contracted by $\scrH(\gamma)_n$ for some $\gamma\in C$.
Recall that~$x$ is contracted if and only if some edge $x\covers y$ is contracted.
Let $x_i>x_{i+1}$, and let~$y$ be obtained from~$x$ by transposing $x_i$ and $x_{i+1}$.
Proposition \ref{associated ji} says that $\scrH(\gamma)_n$ contracts the edge $x\covers y$ if and only if it contracts the join-irreducible 
$\lambda(x,i)$.
By Proposition~\ref{avoid}, $\lambda(x,i)$ is contracted if and only if it contains the pattern $\gamma$.
Since the unique descent in $\lambda(x,i)$ consists of the element $x_i$ followed by the element $x_{i+1}$ and since $\gamma$ is untranslated,
$\lambda(x,i)$ contains $\gamma$ if and only if there is an occurrence of $\gamma$ in $\lambda(x,i)$ which includes the entries $x_i$ and 
$x_{i+1}$.
Also, in the definition of $\lambda(x,i)$, note that all entries of $\lambda(x,i)$ with values weakly between $x_{i+1}$ and $x_i$ are on the same 
side of the pair $(x_i,x_{i+1})$ in~$x$ as in $\lambda(x,i)$.
Thus $\gamma$ occurs in $\lambda(x,i)$ if and only if there is some scramble $\sigma$ of $\lambda$ which occurs in~$x$ such that the cliff
of $\sigma$ occurs in positions $i$ and $i+1$ of~$x$.
Letting $i$ vary over all descents of~$x$, we have that~$x$ is contracted if and only if there is some scramble $\sigma$ of $\lambda$ which 
occurs in~$x$ with adjacent cliff.
\end{proof}

\begin{remark}\rm
\label{monoid}
We now describe how $\K[Z^\Theta_\infty]$ can be obtained via a construction due to Duchamp, Hivert, Novelli and Thibon 
(see~\cite[Proposition 3.12]{NCSymVI}, \cite{arbres} and~\cite[Proposition 18]{Search}).
This construction begins with the free monoid~$M$ on an infinite alphabet and realizes $\K[S_\infty]$ as a Hopf algebra by associating
each permutation~$x$ to the sum of all elements of~$M$ whose ``standardization'' is~$x$.
Given a monoid congruence on~$M$ generated by relations of the form $w\equiv w'$, where $w'$ is obtained from $w$ by transposing 
two adjacent letters, one obtains an equivalence on permutations whenever the congruence on~$M$ is compatible with standardization.
This equivalence defines a sub Hopf algebra of $\K[S_\infty]$ whenever the congruence on~$M$ is compatible with 
``restriction to intervals.''

Starting with an $\scrH$-family $\set{\Theta_n}$ and guided by Theorem~\ref{c-a avoid}, one can construct a congruence on~$M$ which is compatible with standardization and 
restriction to intervals and which recovers the congruences~$\Theta_n$.
Thus the construction via monoid congruences produces a strictly larger class of sub Hopf algebras.
(The example in~\cite[Proposition 3.12]{NCSymVI}, for example, does not correspond to a family of lattice congruences.)
However, there is no immediate way to tell from this construction which sub Hopf algebras can even be described in terms of partial orders,
much less which of them arise from lattice congruences.
Thus, while each $\K[Z^\Theta_\infty]$ arises as a special case of the construction by monoid congruences, 
it is not apparent how one would arrive at the appropriate congruences on~$M$ without the analysis given in the present paper.
\end{remark}

\begin{remark}
\label{compute}
\rm
Computing products in $\K[Z_\infty^\Theta]$ via Equation~(\ref{prod}) involves only identifying permutations with $\pidown x=x$, where $\pidown$ is the downward projection associated to the congruence $\scrH(C)_n$.
This means checking the pattern avoidance condition imposed by Theorem~\ref{c-a avoid}.
However, to compute coproducts by Equation~(\ref{coprod}) one needs to know $\pidown x$ for every $x\in S_n$.
The proof of Theorem~\ref{c-a avoid} indicates how to compute $\pidown x$ inductively.
Suppose for some $\gamma\in C$ with cliff at $j$ that there is a scramble of $\gamma$ occurring as the subsequence 
$x_{i_1}x_{i_2}\cdots x_{i_k}$ of~$x$, with $i_{j+1}=i_j+1$.
Then $x\covers y$ and $x\equiv y$, where~$y$ is obtained from~$x$ by transposing the entries $x_{i_j}$ and $x_{i_j+1}$.
Since $\pidown x=\pidown y$ we continue inductively until we reach an uncontracted permutation.
\end{remark}

\begin{remark}
\label{anti remark}
\rm
The definition of a (not-necessarily graded) Hopf algebra requires the existence of a map $S$ called the {\em antipode}.
However, when $K[X_\infty]$ is a graded, connected Hopf algebra, as defined above, it always possesses an antipode.

Let $\K[X_\infty]$ and $\K[Y_\infty]$ be graded connected Hopf algebras whose antipodes are $S_X$ and $S_Y$.
Suppose $c$ embeds $\K[Y_\infty]$ as a graded sub Hopf algebra of $\K[X_\infty]$.
It is known that $S_Y=r\circ S_X\circ c$, where $r$ is the inverse map, defined on the image of $c$.
In the case of $\K[S_\infty]$ and $\K[Z^\Theta_\infty]$, the maps $c$ and $r$ were defined in Section~\ref{trans sec} and we have 
\begin{equation}
\label{antipode}
S_Z=r\circ S_S\circ c
\end{equation}
In~\cite[Theorem 5.4]{Ag-So}, Aguiar and Sottile give a formula for $S_S$.
This formula and Equation~(\ref{antipode}), along with Remark~\ref{compute}, allows computation of $S_Z$ for any $\Theta$.
\end{remark}

\begin{remark}
\label{combinatorial}
\rm
If $\K[S_\infty]$ is given the structure of a {\em combinatorial} Hopf algebra in the sense of~\cite{CHA}, then this structure can be pulled 
back to $\K[Z^\Theta_\infty]$.
One defines a multiplicative character $\zeta_Z:=\zeta_S\circ c$, where $\zeta_S$ is the chosen multiplicative character of $\K[S_\infty]$, 
so that by definition $c$ is a morphism of combinatorial Hopf algebras.
By \cite[Proposition 5.8(f)]{CHA}, since $c$ is injective, the odd (respectively even) subalgebra of $\K[Z^\Theta_\infty]$ is the image 
under $r$ of the odd (respectively even) subalgebra of $\K[S_\infty]$.
\end{remark}

\section{Examples}
\label{ex}
We conclude by discussing some examples, which by no means exhaust the possibilities.
The examples given in the introduction define $\scrH$-families of congruences.
Specifically, the Tamari lattice is known~\cite{Nonpure II} to be the sublattice of $S_n$ consisting of $312$-avoiding permutations.
The permutation 312 is an untranslated join-irreducible and the only scramble of 312 is 312 itself.
It is easy to check that the pattern $312$ occurs in a permutation~$x$ if and only if it occurs in~$x$ with adjacent cliff.
One can specialize~\cite[Theorems 6.2 and 6.4]{Cambrian} to state that the Tamari lattice is $S_n$ mod the congruence $\scrH(312)_n$, 
or alternately $\scrH(231)_n$.
The fibers of the (left) descent map, in the case of $S_n$, are the congruence classes of $\scrH(\set{231,312})_n$.
Thus we recover the setup described in the introduction.

For a second example, we construct an infinite sequence $\set{\K[S_{\infty,k}]}_{k\ge 1}$ of graded Hopf algebras, each included in its 
successor, limiting to $\K[S_\infty]$, such that the first Hopf algebra consists of one-dimensional graded pieces, and the second Hopf 
algebra is the Hopf algebra of non-commutative symmetric functions.
Since $\scrH_\infty$ is graded by size, for each $k\ge 1$ we define an $\scrH$-family 
$\set{\Delta_{n,k}}_{n\ge 0}:=\set{H(C_k)_n}_{n\ge 0}$ by letting $C_k$ be the set of untranslated join-irreducibles of size $k+1$.
The congruence $\Delta_{n,k}$ is homogeneous of degree $k$, and
by Theorem~\ref{c-a avoid}, $\Delta_{n,k}$ contracts every permutation~$x$ containing adjacent elements $x_i$ and $x_{i+1}$ with 
$x_i-x_{i+1}\ge k$.
Notice that $\Delta_{n,1}$ is the congruence on $S_n$ with a single equivalence class, and $\K[S_{\infty,1}]$ is the binomial Hopf algebra
$B_1$ of~\cite[Section V.2]{Joni-Rota}.
Also, $\Delta_{n,2}$ is the congruence associated to the descent map.
When $k\ge n$, the congruence $\Delta_{n,k}$ consists of congruence classes which are all singletons.

Let $S_{n,k}$ be the quotient of $S_n$ with respect to $\Delta_{n,k}$, and as usual identify $S_{n,k}$ as the induced subposet $\pidown(S_n)$, 
where $\pidown$ is the downward projection corresponding to $\Delta_{n,k}$.
Specifically, $S_{n,k}$ is the subposet of $S_n$ consisting of permutations with no right descents of magnitude $k$ or larger.
Applying Theorem~\ref{subalgebra}, we obtain an infinite sequence $\K[S_{\infty,k}]$ of graded Hopf algebras, each included as a sub 
Hopf algebra of its successor, limiting to the Malvenuto-Reutenauer Hopf algebra $\K[S_\infty]$.

By a simple argument involving induction on $n$, we have
\[|S_{n,k}|=\prod_{i=1}^n\min(i,k).\]
Thinking of $S_n$ as a poset of regions as explained in Section~\ref{symmetric} we see that $\Delta_{n,k}$ contracts every cover 
$R_1\covered R_2$ such that $R_1$ and $R_2$ are separated by a hyperplane normal to $e_a-e_b$ for $a-b\ge k$.
Thus each $S_{n,k}$ is obtained from $S_n$ by deleting hyperplanes, so $S_{n,k}$ is a lattice of regions, or in other words, $S_{n,k}$ 
is obtained by directing the $1$-skeleton of a zonotope.

A related construction yields a sequence of graded Hopf algebras limiting to the Hopf algebra of planar binary trees.
For each $k\ge 1$, let $C'_k$ be the set containing the two permutations $231$ and $(k+1)123\cdots k$ and define an $\scrH$-family 
$\set{\Phi_{n,k}}_{n\ge 0}:=\set{\scrH(C'_k)_n}_{n\ge 0}$.
For each $k$ the family $\set{\Phi_{n,k}}_{n\ge 0}$ is the meet, in the distributive lattice of order ideals of 
$\scrH_\infty$, of $\set{\Delta_{n,k}}_{n\ge 0}$ and $\set{\scrH(231)_n}_{n\ge 0}$.
Construct $P_{n,k}$ and $\K[P_{\infty,k}]$ from $\set{\Phi_{n,k}}$ exactly as $S_{n,k}$ and $\K[S_{\infty,k}]$ were constructed from
$\set{\Delta_{n,k}}$.
We obtain an infinite sequence of Hopf algebras, each included as a sub Hopf algebra of its successor, limiting to the Hopf algebra of planar 
binary trees.
Also, $\K[P_{\infty,k}]$ is a sub Hopf algebra of $\K[S_{\infty,k}]$ for each $k$.
For $k=1,2$ we have $\K[P_{\infty,k}]=\K[S_{\infty,k}]$ and a simple argument shows that the dimensions of the graded pieces of 
$\K[P_{\infty,3}]$ satisfy the defining recurrence of the Pell numbers, sequence A000129 in~\cite{Sloane}.

Our final example is an $\scrH$-family such that the congruence classes appear to be equinumerous with the Baxter permutations.
Say a permutation~$x$ is a {\em twisted Baxter permutation} if and only if the following two conditions hold:
\begin{enumerate}
\item[(i) ] For any 2413-pattern in~$x$, the ``4'' and the ``1'' are not adjacent in~$x$, and
\item[(ii) ] For any 3412-pattern in~$x$, the ``4'' and the ``1'' are not adjacent in~$x$.
\end{enumerate}
For the definition of the usual Baxter permutations, see for example~\cite{CGHK}.
West~\cite{West} showed that the Baxter permutations are the permutations~$x$ such that if 2413 occurs in~$x$ then it occurs as a subpattern 
of some 25314 pattern in~$x$, and if 3142 occurs, then it occurs as a subpattern of some 41352.
This is easily checked to be the set of permutations satisfying (i) above and (ii') below.
\begin{enumerate}
\item[(ii') ] For any 3142-pattern in~$x$, the ``1'' and the ``4'' are not adjacent in~$x$.
\end{enumerate}
Computer calculations show that for $n\le 15$ the twisted Baxter permutations in $S_n$ are equinumerous with the Baxter permutations in 
$S_n$.

The congruences associated to the Tamari lattice and the descent map are both homogeneous of degree two, and these are the only 
$\scrH$-families of homogeneous degree-two congruences.
Specifying an $\scrH$-family of homogeneous degree-three congruences amounts to choosing a non-empty subset of 
$\set{2341,3412,2413,4123}$.
The following theorem is an immediate corollary of Theorem~\ref{c-a avoid}.
\begin{theorem}
\label{2413 3412}
The quotient of $S_n$ mod $\scrH(3412,2413)_n$ is isomorphic to the subposet of $S_n$ induced by the twisted Baxter permutations.
\end{theorem}
Thus by Theorem~\ref{subHopf}, the twisted Baxter permutations are the basis of a Hopf algebra which can be embedded as a sub Hopf 
algebra of the Malvenuto-Reutenauer Hopf algebra.

\begin{prop}
\label{meet}
The congruence $\scrH(3412,2413)_n$ is the meet $\scrH(231)_n\meet\scrH(312)_n$ of the two congruences defining the Tamari lattices.
\end{prop}
The weaker statement, that $\set{\scrH(3412,2413)_n}$ is the meet, among $\scrH$-families, of $\set{\scrH(231)_n}$ and 
$\set{\scrH(312)_n}$ is immediate by inspection of $\scrH_\infty$. 
\begin{proof}
By Theorem~\ref{c-a avoid}, a join-irreducible is contracted by $\scrH(231)_n$ if and only if it contains a 231-pattern, and similarly for 
$\scrH(312)_n$.
Thus the join-irreducibles contracted by $\scrH(231)_n\meet\scrH(312)_n$ are exactly the join-irreducibles containing both a 231- and a 
312-pattern.
Since 2413 and 3412 each contain a 231- and a 312-pattern, by Theorem~\ref{c-a avoid}, every join-irreducible contracted by 
$\scrH(3412,2413)_n$ is also contracted by $\scrH(231)_n\meet\scrH(312)_n$.
Conversely, if $\gamma$ is a join-irreducible contracted by $\scrH(231)_n\meet\scrH(312)_n$, let $A$ be the associated subset.
Since $\gamma$ contains 231, there is an element $c\in A^c$ with $m<c<M$, and since $\gamma$ contains 312, there is an element $b\in A$ 
with $m<b<M$.
If $b<c$ then $cMmb$ is a 3412-pattern in $\gamma$, and if $b>c$ then $cMmb$ is a 2413-pattern in $\gamma$.
Thus $\gamma$ is also contracted by $\scrH(2413,3412)_n$.
\end{proof}

\section{Acknowledgments}
The author wishes to thank Alexander Barvinok, Patricia Hersh, Florent Hivert, Sam Hsiao, Jean-Christophe Novelli, Vic Reiner and 
John Stembridge for helpful conversations.

\newcommand{\journalname}[1]{\textrm{#1}}
\newcommand{\booktitle}[1]{\textrm{#1}}

\end{document}